\documentclass[pdflatex,aap,preprint]{imsart}

\RequirePackage[OT1]{fontenc}
\RequirePackage{amsthm,amsmath}
\RequirePackage[numbers]{natbib}
\RequirePackage[colorlinks,citecolor=blue,urlcolor=blue]{hyperref}
\RequirePackage{hypernat}

\arxiv{math.0905.1863}

\startlocaldefs
\usepackage{latexsym}
\usepackage{graphicx}
\usepackage{dsfont}
\usepackage{amssymb}

\def\be{\begin{eqnarray}}
\def\ee{\end{eqnarray}}

\def\b*{\begin{eqnarray*}}
\def\e*{\end{eqnarray*}}

\newtheorem{Theorem}{Theorem}[section]
\newtheorem{Definition}[Theorem]{Definition}
\newtheorem{Proposition}[Theorem]{Proposition}

\newtheorem{Lemma}[Theorem]{Lemma}
\newtheorem{Corollary}[Theorem]{Corollary}
\newtheorem{Remark}[Theorem]{Remark}
\newtheorem{Example}[Theorem]{Example}

\makeatletter \@addtoreset{equation}{section}

\def \E{\mathbb{E}}
\def \F{\mathbb{F}}

\def \L{\mathbb{L}}
\def \M{\mathbb{M}}
\def \N{\mathbb{N}}
\def \P{\mathbb{P}}

\def \R{\mathbb{R}}
\def \S{\mathbb{S}}

\def\Ac{\mathcal{A}}

\def\Dc{\mathcal{D}}
\def\Ec{\mathcal{E}}
\def\Fc{\mathcal{F}}

\def\Lc{\mathcal{L}}

\def\Rc{\mathcal{R}}
\def\Sc{\mathcal{S}}

\def\vo{\overline{v}}

\def\Xh{{\hat X}}

\def\FF{\mathbf{F}}

\def\TT{\mathbf{T}}

\def\Tr#1{{\rm Tr}\left[#1\right]}

\def \Frac{\displaystyle\frac}

\def \Sup{\displaystyle\sup}

\def\no{\noindent}

\def\={\;=\;}
\def\.{\;.}

\def\eps{\varepsilon}

\def\reff#1{{\rm(\ref{#1})}}

\def\1{\mathds{1}}

\def \ep{\hbox{ }\hfill{ ${\cal t}$~\hspace{-5.1mm}~${\cal u}$   } }

\def \proof{{\noindent \bf Proof. }}
\def \ep{\hbox{ }\hfill$\Box$}

 \def\normeL2#1{\left\|{#1}\right\|_{L^2}}
\newcommand {\lb} {\lambda}
\newcommand {\Chi} {{\bf \raise 1.5pt \hbox{$\delta$}}}

\newcommand{\lap}{\bigtriangleup}
\newcommand{\nb}{\nabla}

\endlocaldefs

\begin{document}

\begin{frontmatter}

\title{A Probabilistic Numerical Method for Fully Nonlinear Parabolic PDEs\protect\thanksref{T1,T2}}
\thankstext{T2}{We are grateful to Mete Soner, Bruno Bouchard, Denis Talay, and Fr\'ed\'eric Bonnans for fruitful comments and suggestions.} 
\thankstext{T1}{This research is part of the Chair {\it Financial Risks} of the {\it Risk Foundation} sponsored by Soci\'et\'e
         G\'en\'erale, the Chair {\it Derivatives of the Future} sponsored by the {F\'ed\'eration Bancaire Fran\c{c}aise}, 
         the Chair {\it Finance and Sustainable Development} sponsored by EDF and Calyon.
}
\runtitle{\tiny{A Probabilistic Numerical Method for Fully Nonlinear Parabolic PDEs}}

\begin{aug}
\author{\fnms{Arash} \snm{Fahim}\thanksref{m1}\ead[label=e1]{arash.fahim@polytechnique.edu}},
\author{\fnms{Nizar} \snm{Touzi}\thanksref{m1}\ead[label=e2]{nizar.touzi@polytechnique.edu}}
\and
\author{\fnms{Xavier} \snm{Warin}\thanksref{m2}
\ead[label=e3]{xavier.warin@edf.fr}}

\runauthor{N. Touzi et al.}

\affiliation{CMAP, Ecole Polytechnique Paris\thanksmark{m1} and EDF R\&D \& FiME, Finance des March\'es de l'Energie Laboratory Paris \ead[label=u1,url]{www.fime-lab.org}\thanksmark{m2}}

\address{CMAP, \'Ecole Polytechnique\\
Route de Saclay\\ 
91128, Palaiseau, France\\
\printead{e1}\\
\phantom{E-mail:\ }\printead*{e2}}

\address{\'Electricit\'e de France - EDF R\&D\\
 L'avenue du G\'en\'eral De Gaulle\\
92141, Clamart, France\\
\printead{e3}\\
\printead{u1}}
\end{aug}

\begin{abstract}
We consider the probabilistic numerical scheme for fully nonlinear PDEs suggested in \cite{cstv}, and show that it can be introduced naturally as a combination of Monte Carlo and finite differences scheme without appealing to the theory of backward stochastic differential equations. Our first main result provides the convergence of the discrete-time approximation and derives a bound on the discretization error in terms of the time step. An explicit implementable scheme requires to approximate the conditional expectation operators involved in the discretization. This induces a further Monte Carlo error. Our second main result is to prove the convergence of the latter approximation scheme, and to derive an upper bound on the approximation error. Numerical experiments are performed for the approximation of the solution of the mean curvature flow equation in dimensions two and three, and for two and five-dimensional (plus time) fully-nonlinear Hamilton-Jacobi-Bellman equations arising in the theory of portfolio optimization in financial mathematics.
\end{abstract}



\end{frontmatter}

\section{Introduction}
\setcounter{equation}{0}
\label{Sec:sectint}

We consider the probabilistic numerical scheme for the approximation of the solution of a fully-nonlinear parabolic Cauchy problem suggested in \cite{cstv}. In the latter paper, a representation of the solution of the PDE is derived in terms of the newly introduced notion of second order backward stochastic differential equations, assuming that the fully-nonlinear parabolic Cauchy problem has a smooth solution. Then, similarly to the case of backward stochastic differential equations which are connected to semi-linear PDEs, this representation suggests a backward probabilistic numerical scheme.

The representation result of \cite{cstv} can be viewed as an extension of the Feynman-Kac representation result, for the linear case, which is widely used in order to approach the numerical approximation problem from the probabilistic viewpoint, and to take advantage of the high dimensional properties of Monte Carlo methods. Previously, the theory of backward stochastic differential equations provided an extension of these approximation methods to the semi-linear case. See for instance Chevance \cite{che}, El Karoui, Peng and Quenez \cite{epq}, Bally and Pag\`es \cite{bal}, Bouchard and Touzi \cite{bt} and Zhang \cite{zhang}. In particular, the latter papers provide the convergence of the ``natural'' discrete-time approximation of the value function and its partial space gradient with the same $L^2$ error of order $\sqrt h$, where $h$ is the length of time step. The discretization involves the computation of conditional expectations, which need to be further approximated in order to result into an implementable scheme. We refer to \cite{bal}, \cite{bt} and \cite{glw} for an complete asymptotic analysis of the approximation, including the regression error.

In this paper, we observe that the backward probabilistic scheme of \cite{cstv} can be introduced naturally without appealing to the notion of backward stochastic differential equation. This is shown is Section \ref{Sec:sectnumericalscheme} where the scheme is decomposed into three steps:
\\
(i) The Monte Carlo step consists in isolating the linear generator of some underlying diffusion process, so as to split the PDE into this linear part and a remaining nonlinear one. 
\\
(ii) Evaluating the PDE along the underlying diffusion process, we obtain a natural discrete-time approximation by using finite differences approximation in the remaining nonlinear part of the equation. 
\\
(iii) Finally, the backward discrete-time approximation obtained by the above steps (i)-(ii) involves the conditional expectation operator which is not computable in explicit form. An implementable probabilistic numerical scheme therefore requires to replace such conditional expectations by a convenient approximation, and induces a further Monte Carlo type of error. 

In the present paper, we do not require the fully nonlinear PDE to have a smooth solution, and we only assume that it satisfies a comparison result in the sense of viscosity solutions. Our main objective is to establish the convergence of this approximation towards the unique viscosity solution of the fully-nonlinear PDE, and to provide an asymptotic analysis of the approximation error.

Our main results are the following. We first prove the convergence of the discrete-time approximation for general nonlinear PDEs, and we provide bounds on the corresponding approximation error for a class of Hamilton-Jacobi-Bellman PDEs. Then, we consider the implementable scheme involving the Monte Carlo error, and we similarly prove a convergence result for general nonlinear PDEs, and we provide bounds on the error of approximation for Hamilton-Jacobi-Bellman PDEs. We observe that our convergence results place some restrictions on the choice of the diffusion of the underlying diffusion process. First, a uniform ellipticity condition is needed; we believe that this technical condition can be relaxed in some future work. More importantly, the diffusion coefficient is needed to dominate the partial gradient of the remaining nonlinearity with respect to its Hessian component. Although we have no theoretical result that this condition is necessary, our numerical experiments show that the violation of this condition leads to a serious mis-performance of the method, see Figure \ref{figfin1}.

Our proofs rely on the monotonic scheme method developed by Barles and Souganidis \cite{barlessouganidis} in the theory of viscosity solutions, and the recent method of shaking coefficients of Krylov \cite{krylov}, \cite{krylov1} and \cite{krylov2} and Barles and Jakobsen \cite{barlesjakobsen}, \cite{barlesjakobsen1} and \cite{barlesjakobsen2}. The use of the latter type of methods in the context of a stochastic scheme seems to be new. Notice however, that our results are of a different nature than the classical error analysis results in the theory of backward stochastic differential equations, as we only study the convergence of the approximation of the value function, and no information is available for its gradient or Hessian with respect to the space variable.

The following are two related numerical methods based on finite differences in the context of Hamilton-Jacobi-Bellman nonlinear PDEs:
\begin{itemize}
\item Bonnans and Zidani \cite{bz} introduced a finite difference scheme which satisfies the crucial monotonicity condition of Barles and Souganidis \cite{barlessouganidis} so as to ensure its convergence. Their main idea is to discretize both time and space, approximate the underlying controlled forward diffusion for each fixed control by a controlled local Markov chain on the grid, approximate the derivatives in certain directions which are found by solving some further optimization problem, and optimize over the control. Beyond the curse of dimensionality problem which is encountered by finite differences schemes, we believe that our method is much simpler as the monotonicity is satisfied without any need to treat separately the linear structures for each fixed control, and without any further investigation of some direction of discretization for the finite differences. 

\item An alternative finite-differences scheme is the semi-Lagrangian method which solves the monotonicity requirement by absorbing the dynamics of the underlying state in the finite difference approximation,  see e.g. Debrabant and Jakobsen \cite{debrabantjakobsen}, Camilli and Jacobsen \cite{camillijakobsen}, Camilli and Falcone \cite{camillifalcone}, and Munos and Zidani \cite{munoszidani} . Loosely speaking, this methods is close in spirit to ours, and corresponds to freezing the Brownian motion $W_h$, over each time step $h$, to its average order $\sqrt{h}$. However it does not involve any simulation technique, and requires the interpolation of the value function at each time step. Thus it is also subject to the curse of dimensionality problems.
\end{itemize}

We finally observe a connection with the recent work of Kohn and Serfaty \cite{ks} who provide a deterministic game theoretic interpretation for fully nonlinear parabolic problems. The game is time limited and consists of two players. At each time step, one tries to maximize her gain and the other to minimize it by imposing a penalty term to her gain. The nonlinearity of the fully nonlinear PDE appears in the penalty. Also, although the nonlinear penalty does not need to be elliptic, a parabolic nonlinearity appears in the limiting PDE. This approach is very similar to the representation of \cite{cstv} where such a parabolic envelope appears in the PDE, and where the Brownian motion plays the role of Nature playing against the player.

The paper is organized as follows. In Section \ref{Sec:sectnumericalscheme}, we provide a natural presentation of the scheme without appealing to the theory of backward stochastic differential equations. Section \ref{Sec:sectasymptotics} is dedicated to the asymptotic analysis of the discrete-time approximation, and contains our first main convergence result and the corresponding error estimate. In Section \ref{Sec:sectprobnum}, we introduce the implementable backward scheme, and we further investigate the induced Monte Carlo error. We again prove convergence and we provide bounds on the approximation error. Finally, Section \ref{Sec:sectnulmerics} contains some numerical results for the mean curvature flow equation on the plane and space, and for a five-dimensional Hamilton-Jacobi-Bellman equation arising in the problem of portfolio optimization in financial mathematics.

\vspace{5mm}

\no {\bf Notations}\quad 
For scalars $a,b\in\R$, we write $a\wedge b:=\min\{a,b\}$, $a\vee b:=\max\{a,b\}$, $a^-:=\max\{-a,0\}$, and $a^+:=\max\{a,0\}$.
By $\M(n,d)$, we denote the collection of all $n\times d$ matrices with real entries. The collection of all symmetric matrices of size $d$ is denoted $\S_d$, and its subset of nonnegative symmetric matrices is denoted by $\S_d^+$.
For a matrix $A\in\M(n,d)$, we denote by $A^{\rm T}$ its transpose. For $A,B\in\M(n,d)$, we denote $A\cdot B:={\rm Tr}[A^{\rm T}B]$. In particular, for $d=1$, $A$ and $B$ are vectors of $\R^n$ and $A\cdot B$ reduces to the Euclidean scalar product.

For a function $u$ from $[0,T]\times\R^d$ to $\R$, we say that $u$ has $q-$polynomial growth (resp. $\alpha-$exponential growth) if
 \b*
 \sup_{t\le T,~x\in\R^d} \frac{|u(t,x)|}{1+|x|^q}<\infty,
 &\mbox{(resp.}&
 \sup_{t\le T,~x\in\R^d} e^{-\alpha |x|}|u(t,x)|<\infty \mbox{)}.
 \e*
For a suitably smooth function $\varphi$ on $Q_T:=(0,T]\times\R^d$, we define 
\b*
|\varphi|_\infty:=\sup_{(t,x)\in Q_T}|\varphi(t,x)|
&\mbox{and}&
|\varphi|_1:=|\varphi|_\infty+\mathop{\sup}\limits_{Q_T\times Q_T}\frac{|\varphi(t,x)-\varphi(t',x')|}{(x-x')+|t-t'|^\frac{1}{2}}.
\e*
Finally, we denote the $\L^p-$norm of a r.v. $R$ by $\Vert R\Vert_p:=\left(\E[|R|^p]\right)^{1/p}$.

\section{Discretization}
\setcounter{equation}{0}
\label{Sec:sectnumericalscheme}

Let $\mu$ and $\sigma$ be two maps from $\R_+\times\R^d$ to $\R^d$ and $\M(d,d)$, respectively. With $a:=\sigma\sigma^{\rm T}$. We define the linear operator:
 \b*
 \Lc^X\varphi
 &:=&
 \frac{\partial\varphi}{\partial t}+\mu\cdot D\varphi+\frac{1}{2}a\cdot D^2\varphi.
 \e*
Given a map 
\b*
F:(t,x,r,p,\gamma)\in\R_+\times\R^d\times\R\times\R^d\times\S_d &\longmapsto& F(x,r,p,\gamma)\in\R
\e*
we consider the Cauchy problem:
\be
\label{equation}
&&-\Lc^X v-F\left(\cdot,v,D v,D^2v\right)= 0,~~\mbox{on}~[0,T)\times\R^d,\\
&&v(T,\cdot)=g,~~\mbox{on}~\in\R^d. \label{terminal}
\ee
Under some conditions, a stochastic representation for the solution of this problem was provided in \cite{cstv} by means of the newly introduced notion of second order backward stochastic differential equations. As an important implication, such a stochastic representation suggests a probabilistic numerical scheme for the above Cauchy problem. 

The chief goal of this section is to obtain the probabilistic numerical scheme suggested in \cite{cstv} by a direct manipulation of \reff{equation}-\reff{terminal} without appealing to the notion of backward stochastic differential equations. 

To do this, we consider an $\R^d$-valued Brownian motion $W$ on a filtered probability space $\left(\Omega,\Fc,\F,\P\right)$, where the filtration $\F=\{\Fc_t,t\in[0,T]\}$ satisfies the usual completeness conditions, and $\Fc_0$ is trivial. 

For a positive integer $n$, let $h:=T/n$, $t_i=ih$, $i=0,\ldots,n$, and consider the one step ahead Euler discretization
 \begin{equation}\label{Euler}
 \hat X_h^{t,x}
 :=
 x+\mu(t,x)h+\sigma(t,x)(W_{t+h}-W_t),
 \end{equation}
of the diffusion $X$ corresponding to the linear operator $\Lc^X$. Our analysis does not require any existence and uniqueness result for the underlying diffusion $X$. However, the subsequent formal discussion assumes it in order to provides a natural justification of our numerical scheme. 

Assuming that the PDE \reff{equation} has a classical solution, it follows from It\^o's formula that
 \b*
 \E_{t_i,x}\left[v\left(t_{i+1},X_{t_{i+1}}\right)\right]
 &=&
 v\left(t_i,x\right)+\E_{t_i,x}\left[\int_{t_i}^{t_{i+1}}\mathcal{L}^X v(t,X_t)dt\right]
 \e*
where we ignored the difficulties related to local martingale part, and $\E_{t_i,x}:=\E[\cdot|X_{t_i}=x]$ denotes the expectation operator conditional on $\{X_{t_i}=x\}$. Since $v$ solves the PDE \reff{equation}, this provides
 \b*
 v(t_i,x)
 &=&
 \E_{t_i,x}\left[v\left(t_{i+1},X_{t_{i+1}}\right)\right]
+\E_{t_i,x}\left[\int_{t_i}^{t_{i+1}}F(\cdot,v,Dv,D^2v)(t,X_t)dt\right].
 \e*
By approximating the Riemann integral, and replacing the process $X$ by its Euler discretization, this suggest the following approximation of the value function $v$
 \be\label{scheme}
 v^h(T,.):=g
 &\mbox{and}&
 v^h(t_i,x):=\TT_h[v^h](t_i,x),
 \ee
where we denoted for a function $\psi:\R_+\times\R^d\longrightarrow\R$ with exponential growth:
 \be\label{TT}
 \TT_h[\psi](t,x)
 :=
 \E\left[\psi(t+h,\hat X_h^{t,x})\right] + hF\left(\cdot,\Dc_h\psi\right)(t,x),
 \\
 \Dc_h^k\psi(t,x)
 :=
 \E[D^k\psi(t+h,\hat X^{t,x}_h)],~~
 k=0,1,2,~~
 \Dc_h\psi
 :=
 \left(\Dc_h^0\psi,\Dc_h^1\psi,\Dc_h^2\psi\right)^\text{T},
 \ee
and $D^k$ is the $k-$th order partial differential operator with respect to the space variable $x$. The differentiations in the above scheme are to be understood in the sense of distributions. This algorithm is well-defined whenever $g$ has exponential growth and $F$ is a Lipschitz map. To see this, observe that any function with exponential growth has weak gradient and Hessian because the Gaussian kernel is a Schwartz function, and the exponential growth is inherited at each time step from the Lipschitz property of $F$. 

At this stage, the above backward algorithm presents the serious drawback of involving the gradient $Dv^h(t_{i+1},.)$ and the Hessian $D^2v^h(t_{i+1},.)$ in order to compute $v^h(t_i,.)$. The following result avoids this difficulty by an easy integration by parts argument.

\begin{Lemma}\label{lemintpart}
For every function $\varphi:Q_T\to\R$ with exponential growth, we have:
 \b*
 \Dc_h\varphi(t_i,x)
 &=&
 \E\left[\varphi(t_{i+1},\hat X^{t_i,x}_h)H_h(t_i,x)\right],
 \e*
where $H_h=(H^h_0,H^h_1,H^h_2)^\text{T}$ and
 \be\label{hermit}
 H^h_0=1,
 &H^h_1=\left({\sigma^{\rm T}}\right)^{-1}\;\Frac{W_h}{h},&
 H^h_2=\left({\sigma^{\rm T}}\right)^{-1}\;\frac{W_hW^{\rm T}_h-h\mathbf{I}_d}{h^2}\;\sigma^{-1}.
 \ee
\end{Lemma}

\proof
The main ingredient is the following easy observation. Let $G$ be a one dimensional Gaussian random variable with unit variance. Then, for any function $f:\R\longrightarrow\R$ with exponential growth, we have:
 \be\label{Hermite0}
 \E[f(G)H^k(G)] &=& \E[f^{(k)}(G)],
 \ee
where $f^{(k)}$ is the $k-$th order derivative of $f$ in the sense of distributions, and $H^k$ is the one-dimensional Hermite polynomial of degree $k$.

\no {\bf 1} \quad Now, let $\varphi:\R^d\longrightarrow\R$ be a function with exponential growth. Then, by direct conditioning, it follows from \reff{Hermite0} that
 \b*
 \E\left[\varphi(\hat X^{t,x}_h)W^i_h\right]
 &=&
 h\sum_{j=1}^d \E\left[\frac{\partial\varphi}{\partial x_j}(\hat X^{t,x}_h)\sigma_{ji}(t,x)\right],
 \e*
and therefore:
 \b*
 \E\left[\varphi(\hat X^{t,x}_h)H^h_1(t,x)\right]
 &=&
 \sigma(t,x)^{\rm T}\E\left[\nb\varphi(\hat X^{t,x}_h)\right].
 \e*
\no {\bf 2}\quad For $i\ne j$, it follows from \reff{Hermite0} that
 \b*
 \E\left[\varphi(\hat X^{t,x}_h)W^i_hW^j_h\right]
 &=&
 h\sum_{k=1}^d \E\left[\frac{\partial\varphi}{\partial x_k}(\hat X^{t,x}_h)W^j_h\sigma_{ki}(t,x)\right]
 \\
 &=&
 h^2\sum_{k,l=1}^d \E\left[\frac{\partial^2\varphi}{\partial x_k\partial x_l}(\hat X^{t,x}_h)\sigma_{lj}(t,x)\sigma_{ki}(t,x)\right],
 \e*
and for $j=i$:
 \b*
 \E\left[\varphi(\hat X^{t,x}_h)((W^i_h)^2-h)\right]
 &=&
 h^2\sum_{k,l=1}^d \E\left[\frac{\partial^2\varphi}{\partial x_k\partial x_l}(\hat X^{t,x}_h)\sigma_{li}(t,x)\sigma_{ki}(t,x)\right].
\e*
This provides:
 \b*
 \E\left[\varphi(\hat X^{t,x}_h)H^h_2(t,x)\right]
 &=&
 \sigma(t,x)^{\rm T}\E\left[\nb^2\varphi(\hat X^{t,x}_h)\sigma(t,x)\right].
 \e*
\ep

\vspace{5mm}

In view of Lemma \ref{lemintpart}, the iteration which computes $v^h(t_i,.)$ out of $v^h(t_{i+1},.)$ in \reff{scheme}-\reff{TT} does not involve the gradient and the Hessian of the latter function.

\begin{Remark}\label{remintpart}{\rm
Clearly, one can proceed to different choices for the integration by parts in Lemma \ref{lemintpart}. One such possibility leads to the representation of $\Dc^h_2\varphi$ as:
 \b*
 \Dc^h_2\varphi(t,x)
 &=&
 \E\left[ \varphi(\hat X^{t,x}_h)(\sigma^{\rm T})^{-1}\frac{W_{h/2}}{(h/2)}
                                 \frac{W^{\rm T}_{h/2}}{(h/2)}\sigma^{-1}
    \right].
 \e* 
This representation shows that the backward scheme \reff{scheme} is very similar to the probabilistic numerical algorithm suggested in \cite{cstv}. 
}
\end{Remark}

Observe that the choice of the drift and the diffusion coefficients $\mu$ and $\sigma$ in the nonlinear PDE \reff{equation} is arbitrary. So far, it has been only used in order to define the underlying diffusion $X$. Our convergence result will however place some restrictions on the choice of the diffusion coefficient, see Remark \ref{remdiffusion}.  

Once the linear operator $\Lc^X$ is chosen in the nonlinear PDE, the above algorithm handles the remaining nonlinearity by the classical finite differences approximation. This connection with finite differences is motivated by the following formal interpretation of Lemma \ref{lemintpart}, where for ease of presentation, we set $d=1$, $\mu\equiv 0$, and $\sigma(x)\equiv 1$:
\begin{itemize}
\item Consider the binomial random walk approximation of the Brownian motion $\hat W_{t_k}:=\sum_{j=1}^k w_j$, $t_k:= kh,$ $k\ge 1$, where $\{w_j, j\ge 1\}$ are independent random variables distributed as $\frac12\left(\delta_{\sqrt{h}}+\delta_{-\sqrt{h}}\right)$. Then, this induces the following approximation:
 \b*
 \Dc^1_h\psi(t,x)
 :=\E\left[\psi(t+h,X^{t,x}_h)H^h_1\right]
 &\approx&
 \frac{\psi(t,x+\sqrt{h})-\psi(t,x-\sqrt{h})}{2\sqrt{h}},
 \e*
which is the centered finite differences approximation of the gradient.
\item Similarly, consider the trinomial random walk approximation $\hat W_{t_k}:=\sum_{j=1}^k w_j$, $t_k:= kh,$ $k\ge 1$, where $\{w_j, j\ge 1\}$ are independent random variables distributed as $\frac{1}{6}\left(\delta_{\{\sqrt{3h}\}}+4\delta_{\{0\}}+\delta_{\{-\sqrt{3h}\}}\right)$ , so that $\E[w_j^n]=\E[W_h^n]$ for all integers $n\le 4$. Then, this induces the following approximation:
 \b*
 \Dc^2_h\psi(t,x)
 &:=&\E\left[\psi(t+h,X^{t,x}_h)H^h_2\right]\\
 &\approx&
 \frac{\psi(t,x+\sqrt{3h})-2\psi(t,x)+\psi(t,x-\sqrt{3h})}{3h},
 \e*
which is the centered finite differences approximation of the Hessian.

\end{itemize}
In view of the above interpretation, the numerical scheme studied in this paper can be viewed as a mixed Monte Carlo--Finite Differences algorithm. The Monte Carlo component of the scheme consists in the choice of an underlying diffusion process $X$. The finite differences component of the scheme consists in approximating the remaining nonlinearity by means of the integration-by-parts formula of Lemma \ref{lemintpart}. 

\section{Asymptotics of the discrete-time approximation}
\setcounter{equation}{0}
\label{Sec:sectasymptotics}

\subsection{The main results}

Our first main convergence results follow the general methodology of Barles and Souganidis \cite{barlessouganidis}, and requires that the nonlinear PDE \reff{equation} satisfies a comparison result in the sense of viscosity solutions. 

We recall that an upper-semicontinuous (resp. lower semicontinuous) function $\underline{v}$ (resp.  $\overline{v}$) on $[0,T]\times\R^d$, is called a viscosity subsolution (resp. supersolution) of \reff{equation} if for any $(t,x)\in[0,T)\times\R^d$ and any smooth function $\varphi$ satisfying
\b*
0=(\underline{v}-\varphi)(t,x)=\!\!\!\!\max_{[0,T]\times\R^d}(\underline{v}-\varphi)\left(\text{resp.}~~0=(\overline{v}-\varphi)(t,x)=\!\!\!\!\min_{[0,T]\times\R^d}(\overline{v}-\psi)\right),
\e*
we have:
\b*
-\Lc^X\varphi-F(t,x,\Dc \varphi(t,x)) 
&\le~\mbox{(resp. $\ge$)}&
0.
\e*

\begin{Definition}\label{defcomp}
We say that \reff{equation} has comparison for bounded functions if for any bounded upper semicontinuous subsolution $\underline{v}$ and any bounded lower semicontinuous supersolution $\overline{v}$ on $[0,T)\times\R^d$, satisfying
\b*
\underline{v}(T,\cdot)\le \overline{v}(T,\cdot),
\e*
we have $\underline{v}\le\overline{v}$.
\end{Definition}

\begin{Remark}\label{rem-comparison}
{\rm Barles and Souganidis \cite{barlessouganidis} use a stronger notion of comparison by accounting for the final condition, thus allowing for a possible boundary layer. In their context, a supersolution $\overline{v}$ and a subsolution $\underline{v}$ satisfy:
\be\label{bscomp1}
\min\left\{-\Lc^X\overline v(T,x)-F(T,x,\Dc\overline v(T,x)) ,\overline v(T,x)-g(x)\right\}
&\le&
0\\
\max\left\{-\Lc^X\underline v(T,x)-F(T,x,\Dc\underline v(T,x)),\underline v(T,x)-g(x)\right\}
&\ge&
0.\label{bscomp2}
\ee
We observe that, by the nature of our equation, \reff{bscomp1} and \reff{bscomp2} imply that the subsolution $\underline{v}\le g$ and the supersolution $\overline{v}\ge g$, i.e. the final condition holds in the usual sense, and no boundary layer can occur. To see this, without loss of generality we suppose that $F(t,x,r,p,\gamma)$ is decreasing with respect to $r$ (see Remark \ref{remFbar}). Let $\varphi$ be a function satisfying 
\b*
0=(\underline{v}-\varphi)(T,x)=\max_{[0,T]\times\R^d}(\underline{v}-\varphi).
\e*
Then define $\varphi_K(t,\cdot)=\varphi(t,\cdot)+K(T-t)$ for $K>0$. Then $\underline{v}-\varphi_K$ also has a maximum at $(T,x)$, and the subsolution property \reff{bscomp1} implies that
\b*
\min\left\{-\Lc^X\varphi(T,x)-F(T,x,\Dc\varphi(T,x))+K ,\underline{v}(T,x)-g(x)\right\}
&\le&
0.
\e*
For a sufficiently large $K$, this provides the required inequality $\underline{v}(T,x)-g(x)\le 0$. A similar argument shows that \reff{bscomp1} implies that $\overline{v}-g\ge 0$. 
}
\end{Remark}

In the sequel, we denote by $F_r$, $F_p$ and $F_\gamma$ the partial gradients of $F$ with respect to $r$, $p$ and $\gamma$, respectively. We also denote by $F_\gamma^-$ the pseudo-inverse of the non-negative symmetric matrix $F_\gamma$. We recall that any Lipschitz function is differentiable a.e. 
 
\vspace{5mm}

\no {\bf Assumption F}\quad {\it {\rm (i)} The nonlinearity $F$ is Lipschitz-continuous with respect to $(x,r,p,\gamma)$ uniformly in $t$, and $|F(\cdot,\cdot,0,0,0)|_\infty<\infty$;
\\
{\rm (ii) $F$ is elliptic and dominated by the diffusion of the linear operator $\Lc^X$, i.e.
 \be\label{F2}
\nb_{\!\!\gamma} F\le a~&\mbox{on}&~\R^d\times\R\times\R^d\times\S_d;
 \ee
{\rm (iii) $F_p\in{\rm Image}(F_\gamma)$ and $\big|F_p^\text{T}F_\gamma^{-}F_p\big|_\infty<+\infty$.
}

\begin{Remark}\label{remF1}{\rm
Assumption {\bf F} (iii) is equivalent to
\be\label{F1}
|m_F^-|_\infty<\infty
&\mbox{where}&
m_F
:=
\min_{w\in\R^d}\left\{F_p\cdot w+w^\text{T}F_\gamma w\right\}.
\ee
This is immediately seen by recalling that, by the symmetric feature of $F_\gamma$, any $w\in\R^d$ has an orthogonal decomposition $w=w_1+w_2\in\text{Ker}(F_\gamma)\oplus\text{Image}(F_\gamma)$, and by the nonnegativity of $F_\gamma$:
 \b*
 F_p\cdot w+w^\text{T}F_\gamma w
 &=&
 F_p\cdot w_1+F_p\cdot w_2+w_2^\text{T}F_\gamma w_2
 \\
 &=&
 -\frac{1}{4}F_p^\text{T}F_\gamma^{-}F_p
 +F_p\cdot w_1
 +\big|\frac12(F_\gamma^{-})^{1/2}\cdot F_p-F_\gamma^{1/2}w_2\big|^2.
 \e*
}
\end{Remark}

\begin{Remark}\label{remdiffusion}{\rm
The above Condition \reff{F2} places some restrictions on the choice of the linear operator $\Lc^X$ in the nonlinear PDE \reff{equation}. First, $F$ is required to be uniformly elliptic, implying an upper bound on the choice of the diffusion matrix $\sigma$. Since $\sigma\sigma^{\rm T}\in\Sc_d^+$, this implies in particular that our main results do not apply to general degenerate nonlinear parabolic PDEs. Second, the diffusion of the linear operator $\sigma$ is required to dominate the nonlinearity $F$ which places implicitly a lower bound on the choice of the diffusion $\sigma$. 
}
\end{Remark}

\begin{Example}{\rm
Let us consider the nonlinear PDE in the one-dimensional case $-\frac{\partial v}{\partial t}-\frac12\left(a^2v_{xx}^+-b^2v_{xx}^-\right)$ where $0<b<a$ are given constants. Then if we restrict the choice of the diffusion to be constant, it follows from Condition {\bf F} that $\frac{1}{3} a^2\le\sigma^2\le b^2$, which implies that $a^2\le 3b^2$. If the parameters $a$ and $b$ do not satisfy the latter condition, then the diffusion $\sigma$ has to be chosen to be state and time dependent.
}
\end{Example}

\begin{Theorem}[Convergence]\label{thmconv}
Let Assumption {\bf F} hold true, and $|\mu|_1$, $|\sigma|_1<\infty$ and $\sigma$ is invertible. Also assume that the fully nonlinear PDE \reff{equation} has comparison for bounded functions. Then for every bounded Lipschitz function $g$, there exists a bounded function $v$ so that
 \b*
 v^h \longrightarrow v
 &&
 \mbox{locally uniformly}.
 \e*
In addition, $v$ is the unique bounded viscosity solution of problem \reff{equation}-\reff{terminal}.
\end{Theorem}

\begin{Remark}
\label{general}{\rm Under the boundedness condition on the coefficients $\mu$ and $\sigma$, the restriction to a bounded terminal data $g$ in the above Theorem \ref{thmconv} can be relaxed by an immediate change of variable. Let $g$ be a function with $\alpha-$exponential growth for some $\alpha>0$. Fix some $M>0$, and let $\rho$ be an arbitrary smooth positive function with:
 \b*
 \rho(x)=e^{\alpha|x|} &\mbox{for}& |x|\ge M,
 \e*
so that both $\rho(x)^{-1}\nb\rho(x)$ and $\rho(x)^{-1}\nb^2\rho(x)$ are bounded. Let 
 \b*
 u(t,x) \;:=\; \rho(x)^{-1}v(t,x)
 &\mbox{for}&
 (t,x)\in[0,T]\times\R^d.
 \e*
Then, the nonlinear PDE problem \reff{equation}-\reff{terminal} satisfied by $v$ converts into the following nonlinear PDE for $u$:
\be \label{pde1}
-\Lc^X u-\tilde F\left(\cdot,u,Du,D^2u\right)=0&&~\mbox{on}~[0,T)\times\R^d\\ 
v(T,\cdot)=\tilde g:=\rho^{-1}g&&~\mbox{on}~\R^d,\nonumber
\ee
where 
\b*
\tilde F(t,x,r,p,\gamma)
&:=&
r\mu(x)\cdot\rho^{-1}\nabla\rho +\frac12{\rm Tr}\left[a(x)\left(r\rho^{-1}\nabla^2\rho +2p \rho^{-1}\nabla\rho^{\rm T}\right)\right] \\
&&
+\rho^{-1} F \left(t,x,r\rho ,r\nabla\rho+p\rho, r\nabla^2\rho+2p\nabla\rho^{\rm T}+\rho\gamma\right) .
\e*
Recall that the coefficients $\mu$ and $\sigma$ are assumed to be bounded. Then, it is easy to see that $\tilde F$ satisfies the same conditions as $F$. Since $\tilde g$ is bounded, the convergence Theorem \ref{thmconv} applies to the nonlinear PDE \reff{pde1}.
\ep}
\end{Remark}

\begin{Remark}{\rm
Theorem \ref{thmconv} states that the inequality \reff{F2} (i.e. diffusion must dominate the nonlinearity in $\gamma$) is sufficient for the convergence of the Monte Carlo--Finite Differences scheme. We do not know whether this condition is necessary:
\\
$\bullet$ Subsection \ref{subsectlinear} suggests that this condition is not sharp in the simple linear case,
\\
$\bullet$ however, our numerical experiments of Section \ref{Sec:sectnulmerics} reveal that the method may have a poor performance in the absence of this condition, see Figure \ref{figfin1}.
}
\end{Remark}

We next provide bounds on the rate of convergence of the Monte Carlo--Finite Differences scheme in the context of nonlinear PDEs of the Hamilton-Jacobi-Bellman type in the same context as \cite{barlesjakobsen}. The following assumptions are stronger than Assumption $\mathbf{F}$ and imply that the nonlinear PDE \reff{equation} satisfies a comparison result for bounded functions.

\vspace{5mm}

\no {\bf Assumption HJB}\quad {\it The nonlinearity $F$ satisfies Assumption {\rm\bf F}(ii)-(iii), and is of the Hamilton-Jacobi-Bellman type:
 \begin{align*}
\frac12 a\cdot\gamma&+b\cdot p+ F(t,x,r,p,\gamma)
 =
 \inf_{\alpha\in\mathcal{A} }\{\mathcal{L}^{\alpha}(t,x,r,p,\gamma)\}
 \\
 \mathcal{L}^{\alpha}(t,x,r,p,\gamma)
 &:=
 \frac{1}{2}
 Tr[\sigma^\alpha\sigma^{\alpha{\rm T}}(t,x)\gamma]
 +b^{\alpha}(t,x)p+c^{\alpha}(t,x)r+f^{\alpha}(t,x)
 \end{align*}
where the functions $\mu$, $\sigma$, $\sigma^\alpha$, $b^\alpha$, $c^\alpha$ and $f^\alpha$ satisfy:
 \b*
 |\mu|_\infty+|\sigma|_\infty+\sup_{\alpha\in\Ac}\left(|\sigma^\alpha|_1+|b^\alpha|_1+|c^\alpha|_1+|f^\alpha|_1\right) &<& \infty.
 \e*}
 
\hspace{5mm}

\no {\bf Assumption HJB+}\quad {\it The nonlinearity $F$ satisfies {\rm \bf HJB}, and for any $\delta>0$, there exists a finite set $\left\{\alpha_i\right\}_{i=1}^{M_\delta}$ such that for any $\alpha\in\Ac$:
 \b*
 \mathop{\inf}\limits_{1\le i\le M_\delta}
 |\sigma^\alpha-\sigma^{\alpha_i}|_\infty+|b^\alpha-b^{\alpha_i}|_\infty+|c^\alpha-c^{\alpha_i}|_\infty+|f^\alpha-f^{\alpha_i}|_\infty
 &\le& \delta.
 \e*
}

\begin{Remark}{\rm
The assumption HJB+ is satisfied if $\Ac$ is a separable topological space and $\sigma^\alpha(\cdot)$, $b^\alpha(\cdot)$, $c^\alpha(\cdot)$ and $f^\alpha(\cdot)$  are continuous maps from $\Ac$ to $C_b^{\frac{1}{2},1}$; the space of bounded maps which are Lipschitz in $x$ and $\frac{1}{2}$--H\"older in $t$.
}
\end{Remark}

\begin{Theorem}[Rate of Convergence]\label{thmrateconv}
Assume that the final condition $g$ is bounded Lipschitz-continuous. Then, there is a constant $C>0$ such that:
\begin{enumerate} 
\item[{\rm(i)}] under Assumption {\bf\rm HJB}, we have $v-v^h \le Ch^{1/4}$,
\item[{\rm(ii)}] under the stronger condition {\bf\rm HJB+}, we have $-Ch^{1/10}\le v-v^h \le Ch^{1/4}$.
\end{enumerate}
\end{Theorem}

The above bounds can be improved in some specific examples. See Subsection \ref{subsectlinear} for the linear case where the rate of convergence is improved to $\sqrt{h}$.

We also observe that, in the PDE Finite Differences literature, the rate of convergence is usually stated in terms of the discretization in the space variable $|\Delta x|$. In our context of stochastic differential equation, notice that $|\Delta x|$ is or the order of $h^{1/2}$. Therefore, the above upper and lower bounds on the rate of convergence corresponds to the classical rate $|\Delta x|^{1/2}$ and $|\Delta x|^{1/5}$, respectively.

\subsection{Proof of the convergence result}

We now provide the proof Theorem \ref{thmconv} by building on Theorem 2.1 and Remark 2.1 of Barles and Souganidis \cite{barlessouganidis} which requires the scheme to be consistent, monotone and stable. Moreover, since we are assuming the (weak) comparison for the equation, we also need to prove that our scheme produces a limit which satisfies the terminal condition in the usual sense, see Remark \ref{rem-comparison}.

Throughout this section, all the conditions of Theorem \ref{thmconv} are in force.

\begin{Lemma}\label{lemcons}
Let $\varphi$ be a smooth function with bounded derivatives. Then for all $(t,x)\in[0,T]\times\R^d$:
 \b*
 \lim_{\tiny{\begin{array}{c}
             (t',x')\to(t,x)
             \\
             (h,c)\to (0,0)
             \\
             t'+h\le T
             \end{array}}}
 \!\!\!\!\!\!\!\!\!\!\!\!\frac{[c+\varphi](t',x')\!-\!\TT_h[c+\varphi](t',x')}{h}
 \!\!\!\!&=&\!\!\!\!
 -\left(\mathcal{L}^X\varphi\!+\!F(\cdot,\varphi,D\varphi,D^2\varphi)\right)(t,x).
 \e*
\end{Lemma}

The proof is a straightforward application of It\^o's formula, and is omitted.

\vspace{5mm}

\begin{Lemma}\label{lemmon}
Let $\varphi, \psi:~[0,T]\times\R^d\longrightarrow\R$ be two Lipschitz functions. Then, for some $C>0$:
 \b*
 \varphi\le\psi
 &\Longrightarrow&
 \TT_h[\varphi](t,x)\le \TT_h[\psi](t,x)+Ch\;\E[(\psi-\varphi)(t+h,\hat X^{t,x}_h)]
  \e*
where $C$ depends only on constant $K$ in \reff{F1}.
\end{Lemma}
\proof
By Lemma \ref{lemintpart} the operator $\TT_h$ can be written as:
 \b*
 \TT_h[\psi](t,x)
 &=&
 \E\left[\psi(\hat X_h^{t,x})\right]+hF\left(t,x,\E[\psi(\hat X_h^{t,x})H_h(t,x)]\right).
 \e*
Let $f:=\psi-\varphi\ge 0$ where $\varphi$ and $\psi$ are as in the statement of the lemma. Let $F_\tau$ denote the partial gradient with respect to $\tau=(r,p,\gamma)$. By the mean value Theorem:
\b*
\TT_h[\psi](t,x)-\TT_h[\varphi](t,x)
&=&
\E\left[f(\hat X_h^{t,x})\right]+h F_\tau(\theta)\cdot\Dc_hf(\hat X_h^{t,x})
\\
&=&
\E\left[f(\hat X_h^{t,x})\left(1+h F_\tau(\theta)\cdot H_h(t,x)\right)\right],
\e*
for some $\theta=(t,x,\bar r,\bar p,\bar \gamma)$. By the definition of $H_h(t,x)$:
\begin{align*}
\TT_h[\psi]-\TT_h[\varphi]
=&
\E\Bigl[f(\hat X_h^{t,x})\Bigl(1+hF_r+F_p.(\sigma^{\rm T})^{-1}W_h\\
     &                       +h^{-1}F_\gamma\cdot(\sigma^{\rm T})^{-1}(W_hW^{\rm T}_h-hI)\sigma^{-1}
                     \Bigr)\Bigr],
\end{align*}
where the dependence on $\theta$ and $x$ has been omitted for notational simplicity. Since $F_\gamma\le a$ by \reff{F1} of Assumption $\FF$, we have $1-a^{-1}\cdot F_\gamma\ge 0$ and therefore:
 \begin{align*}
 \TT_h[\psi]\!\!-\!\!\TT_h[\varphi]\!\!
 \ge&
 \E\Bigl[\!f(\hat X_h^{t,x})\Bigl(h F_r
                         \!   +\!\!F_p\cdot({\sigma^{\rm T}})^{-1}W_h
\!     +\!\!h^{-1} F_\gamma\cdot({\sigma^{\rm T}})^{-1}W_hW^{\rm T}_h\sigma^{-1}
                      \Bigr)\!\Bigr]
 \\
 =&
 \E\Bigl[\!f(\hat X_h^{t,x})\Bigl(h F_r
                           \! +\!\!hF_p\cdot({\sigma^{\rm T}})^{-1}\frac{W_h}{h}
                          \!+\!\!hF_\gamma\cdot({\sigma^{\rm T}})^{-1}\frac{W_hW^{\rm T}_h}{h^2}\sigma^{-1}\Bigr)
   \!\Bigr].
 \end{align*}
Let $m_F^-:=\max\{-m_F,0\}$, where the function $m_F$ is defined in \reff{F1}. Under Assumption {\bf F}, we have $K:=|m_F^-|_\infty<\infty$, then
\b*
F_p.{\sigma^{\rm T}}^{-1}\frac{W_h}{h}
                            +hF_\gamma\cdot{\sigma^{\rm T}}^{-1}\frac{W_hW^{\rm T}_h}{h^2}\sigma^{-1}\ge -K
\e*
one can write,
 \b*
 \TT_h[\psi]-\TT_h[\varphi]
 &\ge&
 \E\left[f(\hat X_h^{t,x})\left(hF_r-hK\right)\right]
 \;\ge\; 
 -C'h\E\left[f(\hat X_h^{t,x})\right]
 \e*
for some constant $C>0$, where the last inequality follows from \reff{F1}.
\ep

\vspace*{5mm}
The following observation will be used in the proof of Theorem \ref{thmrateconv} below.

\begin{Remark}\label{remFbar}
{\rm
The monotonicity result of the previous Lemma \ref{lemmon} is slightly different from that required in \cite{barlessouganidis}. However, as it is observed in Remark 2.1 in \cite{barlessouganidis}, their convergence theorem holds under this approximate monotonicity. From the previous proof, we observe that if the function $F$ satisfies the condition:
 \be\label{strict}
 F_r-\Frac{1}{4}F_p^{\rm T}F_\gamma^{-}F_p
 &\ge&
 0,
 \ee
then, the standard monotonicity condition
 \be\label{standardmonotonicity}
 \varphi \;\le\; \psi
 &\Longrightarrow& 
 \TT_h[\varphi](t,x) \;\le\; \TT_h[\psi](t,x)
 \ee
holds. Using the parabolic feature of the equation, we may introduce a new function $u(t,x):=e^{\theta(T-t)}v(t,x)$ which solves a nonlinear PDE satisfying \reff{strict}. Indeed, direct calculation shows that the PDE inherited by $u$ is:
 \be
 \label{mani}
 -\mathcal{L}^X u-\overline{F}\left(\cdot,u,Du,D^2u\right)= 0,
 &\mbox{on}& 
 [0,T)\times\R^d\\
 u(T,x)=g(x),
 &\mbox{on}& 
 \R^d,\label{termani}
 \ee 
where $\overline{F}(t,x,r,p,\gamma)=e^{\theta(T-t)}F(t,x,e^{-\theta(T-t)}r,e^{-\theta(T-t)}p,e^{-\theta(T-t)}\gamma)+\theta r$. Then, it is easily seen that $\overline{F}$ satisfies the same conditions as $F$ together with \reff{strict} for sufficiently large $\theta$.
}
\end{Remark}

\begin{Lemma}\label{lemstab}
Let $\varphi,\psi:[0,T]\times\R^d\longrightarrow\R$ be two $\L^\infty-$bounded functions. Then there exists a constant $C>0$ such that
 \b*
|{\bf T}_h[\varphi]-{\bf T}_h[\psi]|_\infty\le|\varphi-\psi|_\infty(1+Ch)
 \e*
In particular, if $g$ is $\L^\infty-$bounded, the family $(v^h)_h$ defined in \reff{scheme} is $\L^\infty-$bounded, uniformly in $h$.
\end{Lemma}

\proof
Let $f:=\varphi-\psi$. Then, arguing as in the previous proof, 
\b*
{\bf T}_h[\varphi]-{\bf T}_h[\psi]
\!\!&\!\!=\!\!&\!\!
\E\left[ f(\hat X_h)\left(1-a^{-1}\cdot F_\gamma+h|A_h|^2
                                    +h F_r-\frac{h}{4}F_p^{\rm T}F_\gamma^{-}F_p\right)\right].
\e*
where 
\b*
A_h=\frac12(F_\gamma^{-})^{1/2} F_p-F_\gamma^{1/2}{\sigma^{\rm T}}^{-1}\frac{W_h}{h}.
\e*
Since $1-{\rm Tr}[a^{-1}F_\gamma]\ge 0$, $|F_r|_\infty<\infty$, and $|F_p^{\rm T}F_\gamma^{-}F_p |_\infty<\infty$ by Assumption $\FF$, it follows that
 \b*
 \left|\TT_h[\varphi]-\TT_h[\psi]\right|_\infty
 &\le& 
 | f|_\infty \left( 1-a^{-1}\cdot F_\gamma+h\E[|A_h|^2]+Ch \right)
\e*
But, $\E[|A_h|^2]=\frac{h}{4}F_p^{\rm T}F_\gamma^{-}F_p+a^{-1}\cdot F_\gamma$. Therefore, by Assumption {\bf F}
\b*
 \left|\TT_h[\varphi]-\TT_h[\psi]\right|_\infty
  &\le&
 | f|_\infty\left( 1+\frac{h}{4}F_p^{\rm T}F_\gamma^{-}F_p +Ch\right)
 \;\le\;
 | f|_\infty(1+\bar Ch).
 \e*
To prove that the family $(v^h)_h$ is bounded, we proceed by backward induction. By the assumption of the lemma $v^h(T,.)=g$ is $\L^\infty-$bounded. We next fix some $i<n$ and we assume that $|v^h(t_j,.)|_\infty\le C_j$ for every $i+1\le j\le n-1$. Proceeding as in the proof of Lemma \ref{lemmon} with $\varphi\equiv v^h(t_{i+1},.)$ and $\psi\equiv 0$, we see that
\b*
\left|v^h(t_i,.)\right|_\infty
&\le&
h\left|F(t,x,0,0,0)\right|
+
C_{i+1}(1+Ch).
 \e*
Since $F(t,x,0,0,0)$ is bounded by Assumption $\FF$, it follows from the discrete Gronwall inequality that $|v^h(t_i,.)|_\infty
\le Ce^{CT}$
for some constant $C$ independent of $h$.
\ep

\begin{Remark}\label{reminter}{\rm
The approximate function $v^h$ defined by \reff{scheme} is only defined on $\{ih|i=0,\cdots,N\}\times\R^d$. Our methodology requires to extend it to any $t\in[0,T]$. This can be achieved by any interpolation, as long as the regularity property of $v^h$ mentioned in Lemma \ref{lemregvh} below is preserved. For instance, on may simply use linear interpolation.}
\end{Remark}

\begin{Lemma}\label{lemregvh}
The function $v^h$ is Lipschitz in $x$, uniformly in $h$.
\end{Lemma}

\proof We report the following calculation in the one-dimensional case $d=1$ in order to simplify the presentation. 
\\
{\bf 1.} For fixed $t\in[0,T-h]$, we argue as in the proof of Lemma \ref{lemmon} to see that for $x,x'\in\R^d$ with $x>x'$:
\be\label{lipxexpand}
 v^h(t,x)- v^h(t,x')
 &=&
 A+hB,
 \ee
where, denoting $\delta^{(k)}:=D^kv^h(t+h,\hat X^{t,x}_h)-D^kv^h(t+h,\hat X^{t,x'}_h)$ for $k=0,1,2$:
 \b*
 A
\!\! &\!\!:=\!\!&\!\!
 \E\big[\delta^{(0)}\big]
 +h\biggl( F\Bigl(t,x',\Dc v^h(t+h,\Xh_h^{t,x})\Bigr) \!\!-\!\! F\Bigl(t,x',\Dc v^h(t+h,\Xh_h^{t,x'})\Big)\biggr)
 \\
\!\! &\!\!=\!\!&\!\!
 \E\big[(1+hF_r)\delta^{(0)}+hF_p \delta^{(1)}+hF_\gamma\delta^{(2)}\big],
 \\
 |B|
\!\! &\!\!:=\!\!&\!\!
 \Big|F\Bigl(t,x,\Dc v^h(t+h,\Xh_h^{t,x})\Bigr) \!\!-\!\! F\Bigl(t,x',\Dc v^h(t+h,\Xh_h^{t,x})\Bigr)\Big|
\!\! \le\!\!
 |F_x|_\infty|x-x'|,
 \e*
by Assumption {\bf F} (i). By Lemma \ref{lemintpart} we write for $k=1,2$:
\begin{align*}
\E\big[\delta^{(k)}\big]
\!&\!=\!\!
\E\big[\delta^{(0)}H_k^h(t,x)+v^h(t+h,\hat X^{t,x'}_h)\left(H^h_k(t,x)-H^h_k(t,x')\right)\big]
\\
\!&\!\!=\!\!
\E\big[\delta^{(0)}H_k^h(t,x)
       +Dv^h(t+h,\hat X^{t,x'}_h)\!\!\left(\frac{W_h}{h}\right)^{\!\!k-1}\!\!\!\!\!\!\!\!\!\bigl(\sigma(t,x)^{\!-k}\!\!-\!\!\sigma(t,x')^{\!-k}\bigr)\sigma(t,x')\big].
\end{align*}
Then, dividing both sides of \reff{lipxexpand} by $x-x'$ and taking limsup, if follows from the above equalities that
\b*
\lefteqn{\limsup_{|x-x'|\searrow0}\frac{ |v^h(t,x)- v^h(t,x')|}{(x-x')}}\\
\!\!\!\!&\!\!\le\!\!&\!\!
\E\biggl[\biggl| \limsup_{|x-x'|\searrow0} \frac{v^h(t\!\!+\!\!h,\hat X^{t,x}_h)\!\!-\!\! v^h(t\!\!+\!\!h,\hat X^{t,x'}_h)}{(x-x')}\Bigl(1\!\!+\!\!hF_r\!\!+\!\!F_p\frac{W_h}{\sigma(t,x)}\!\!+\!\!F_\gamma\frac{W_h^2-h}{\sigma(t,x)^2h}\Bigr)\\
&&\hspace{8mm}+Dv^h(t+h,\hat X^{t,x}_h)
               \left(W_hF_\gamma\frac{-2\sigma_x(t,x)}{\sigma(t,x)^{2}}
                     +hF_p\frac{\sigma_x(t,x)}{\sigma(t,x)}
               \right)\biggr|\biggr]+Ch.
\e*
{\bf 2.} Assume $v^h(t+h,.)$ is Lipschitz with constant $L_{t+h}$. Then
\b*
\lefteqn{\limsup_{|x-x'|\searrow0}\frac{ |v^h(t,x)- v^h(t,x')|}{(x-x')}}\\
\!\!&\!\!\le\!\!&\!\!
L_{t+h}\E\biggl[\biggl|(1\!\!+\!\!\mu_x(t,x)h\!\!+\!\!\sigma_x(t,x)\sqrt{h}N)\Bigl(1\!\!+\!\!hF_r\!\!+\!\!\frac{F_p\sqrt{h}N}{\sigma(t,x)}\!\!+\!\!\frac{F_\gamma N^2}{\sigma(t,x)^2}\!\!-\!\!\frac{F_\gamma}{\sigma(t,x)^2}\Bigr)\\
&&\hspace{20mm}+\sqrt{h}NF_\gamma\frac{-2\sigma_x(t,x)}{\sigma(t,x)^{2}}+hF_p\frac{\sigma_x(t,x)}{\sigma(t,x)}\biggr|\biggr]+Ch.
\e*
Observe that 
$$
F_p\frac{\sigma_x}{\sigma}
= 
\sigma_x \frac{F_p}{\sqrt{F_\gamma}} \frac{\sqrt{F_\gamma}}{\sigma}\;\1_{F_\gamma\neq 0}.
$$
Since all terms on the right hand-side are bounded, under our assumptions, it follows that $|F_p\frac{\sigma_x}{\sigma}|_\infty<\infty$ (we emphasize that the geometric structure imposed in Assumption {\bf F} (iii) provides this result in any dimension). Then:
\b*
\lefteqn{\limsup_{|x-x'|\searrow0}\frac{ |v^h(t,x)- v^h(t,x')|}{(x-x')}}\\
\!\!&\!\!\le\!\!&\!\!
L_{t+h}\biggl(\E\biggl[\Bigl|(1\!\!+\!\!\mu_x(t,x)h\!\!+\!\!\sigma_x(t,x)\sqrt{h}N)\Bigl(1\!\!+\!\!\frac{F_p\sqrt{h}N}{\sigma(t,x)}\!\!+\!\!\frac{F_\gamma N^2}{\sigma(t,x)^2}\!\!-\!\!\frac{F_\gamma}{\sigma(t,x)^2}\Bigr)\\
&&\hspace{20mm}+\sqrt{h}NF_\gamma\frac{-2\sigma_x(t,x)}{\sigma(t,x)^{2}}\Bigr|\biggr]+Ch\biggr)+Ch.
\e*
{\bf 3.} Let $\tilde{\P}$ be the probability measure equivalent to $\P$ defined by the density 
 \b*
 Z:=1-\alpha + \alpha N^2
 &\mbox{where}&
 \alpha=\frac{F_\gamma}{\sigma(t,x)^2}.
 \e*
Then,
\begin{align*}
\limsup_{|x-x'|\searrow0}&\frac{ |v^h(t,x)- v^h(t,x')|}{(x-x')}
\!\le\!
L_{t+h}\biggl(\E^{\tilde\P}\biggl[\biggl|\Bigl(1+\mu_x(t,x)h+\sigma_x(t,x)\sqrt{h}N\Bigr)\\
\!\!&\times\Bigl(1+Z^{-1}F_p\frac{\sqrt{h}N}{\sigma(t,x)}\Bigr)+Z^{-1}\sqrt{h}NF_\gamma\frac{-2\sigma_x(t,x)}{\sigma(t,x)^2}\biggr|\biggr]\!+\!Ch\biggr)\!+\!Ch.
\end{align*}
By Cauchy--Schwartz inequality, we have
\begin{align*}
\limsup_{|x-x'|\searrow0}&\frac{ |v^h(t,x)- v^h(t,x')|}{x-x'}
\!\!\le\!\
L_{t+h}\biggl(\E^{\tilde\P}\biggl[\biggl|\Bigl(1+\mu_x(t,x)h+\sigma_x(t,x)\sqrt{h}N\Bigr)\\
\!\!&\times\Bigl(1+Z^{-1}F_p\frac{\sqrt{h}N}{\sigma(t,x)}\Bigr)+Z^{-1}\sqrt{h}NF_\gamma\frac{-2\sigma_x(t,x)}{\sigma(t,x)^{2}}\biggr|^2\biggr]^\frac12\!+\!Ch\biggr)\!+\!Ch.
\end{align*}
By writing back the expectation in terms of probability $\P$,
\b*
\limsup_{|x-x'|\searrow0}&\frac{ |v^h(t,x)- v^h(t,x')|}{x-x'}
\!\!\le\!\!
L_{t+h}\biggl(\E\biggl[Z\biggl|\Bigl(1+\mu_x(t,x)h+\sigma_x(t,x)\sqrt{h}N\Bigr)\\
\!\!\!&\times\Bigl(1+Z^{-1}F_p\frac{\sqrt{h}N}{\sigma(t,x)}\Bigr)+Z^{-1}\sqrt{h}NF_\gamma\frac{-2\sigma_x(t,x)}{\sigma(t,x)^{2}}\biggr|^2\biggr]^\frac12+Ch\biggr)+Ch.
\e*
By expanding the quadratic term inside the expectation, we observe that expectation of all the terms having $\sqrt{h}$, is zero. Therefore,
\b*
\limsup_{|x-x'|\searrow0}&\frac{ |v^h(t,x)- v^h(t,x')|}{(x-x')}
\!\!\le\!\!
L_{t+h}\biggl(\E^{\tilde\P}\biggl[\biggl|\Bigl(1+\mu_x(t,x)h+\sigma_x(t,x)\sqrt{h}N\Bigr)\\
\!\!\!&\times\Bigl(1+Z^{-1}F_p\frac{\sqrt{h}N}{\sigma(t,x)}\Bigr)+Z^{-1}\sqrt{h}NF_\gamma\frac{-2\sigma_x(t,x)}{\sigma(t,x)^{2}}\biggr|^2\biggr]^\frac12+Ch\biggr)+Ch\\
&\le\!\!
L_{t+h}\Bigl((1+C'h)^\frac12+Ch\Bigr)+Ch,
\e*
which leads to 
\b*
\limsup_{|x-x'|\searrow0}\frac{ |v^h(t,x)- v^h(t,x')|}{(x-x')}
&\le&
Ce^{C'T/2},
\e*
for some constants $C,C'>0$.\ep\\

Finally, we prove that the terminal condition is preserved by our scheme as the time step shrinks to zero.

\begin{Lemma}\label{lemsqrtT-t}
For each $x\in\R^d$ and $t_k=kh$ with $k=1,\cdots,n$, we have;
\b*
|v^h(t_k,x)-g(x)|
&\le&
C(T-t_k)^{\frac12}.
\e*
\end{Lemma}
\proof
{\bf 1.} By the same argument as in the proof of Lemma \ref{lemstab}, we have:
and for $j\ge i$:
\b*
v^h(t_{j},\Xh_{t_{j}}^{t_i,x})
\!\!&\!\!=\!\!&\!\!
\E_{t_j}\left[v^h(t_{j+1},\Xh_{t_{j+1}}^{t_i,x})\left(1-\alpha_j+\alpha_jN_j^2\right)\right]\\
\!\!&\!\!\!\!& \!\!+h\biggl(F^j_0\!+\!F^j_r\E_{t_j}[v^h(t_{j+1},\Xh_{t_{j+1}}^{t_i,x})]\!+\!F^j_p\cdot \E_{t_j}[Dv^h(t_{j+1},\Xh_{t_{j+1}}^{t_i,x})]\biggr),
\e*
where $F^j_0:=F(t_{j},\Xh_{t_{j}}^{t_i,x},0,0,0)$, $\alpha_j$, $F^j_r$, $F^j_p$ are $\Fc_{t_j}-$adapted random variables defined as in the proof of Lemma \ref{lemstab} at $t_j$, and $N_j=\frac{W_{t_{j+1}}-W_{t_{j}}}{\sqrt{h}}$ has a standard Gaussian distribution. Combine the above formula for $j$ from $i$ to $n-1$, we see that
 \begin{align*}
 v^h(t_i,x)
 \E\left[g(\Xh_T^{t_i,x})P_{i,n}\right]
 = 
 +h\E\sum_{j=i}^{n-1}F^j_0&+ F^j_r\E_{t_j}[v^h(t_{j+1},\Xh_{t_{j+1}}^{t_i,x})]\\
                          &+F^j_p\cdot \E_{t_j}[Dv^h(t_{j+1},\Xh_{t_{j+1}}^{t_i,x})],
 \end{align*}
where $P_{i,k}:=\prod_{j=i}^{k-1}\left(1-\alpha_j+\alpha_jN_j^2\right)> 0$ a.s. for all $1\le i< k\le n$ and $P_{i,i}=1$. Obviously $\{P_{i,k},i\le k\le n\}$ is a martingale for all $i\le n$, a property which will be used later. Since $|F(\cdot,\cdot,0,0,0)|_\infty<+\infty$, and using Assumption {\bf F} and Lemmas \ref{lemregvh} and \ref{lemstab}:
\be\label{inequality}
|v^h(t_i,x)-g(x)|
&\le&
\left |\E\left[\left(g(\Xh_T^{t_i,x})-g(x)\right)P_{i,n}\right]\right |+C(T-t_i).
\ee
{\bf 2.} Let $\{g_\eps\}_\eps$ be the family of smooth functions obtained from $g$ by convolution with a family of mollifiers $\{\rho_\eps\}$, i.e. $g_\eps=g\ast\rho_\eps$. Note that we have 
 \be\label{geps-g}
 |g_\eps-g|_\infty\le C\eps,~~
 |Dg_\eps|_\infty\le|Dg|_\infty
 &\mbox{and}&
 |D^2g_\eps|_\infty\le\eps^{-1}|Dg|_\infty.
 \ee
Then:
 \be
 \left |\E\left[\left(g(\Xh_T^{t_i,x})-g(x)\right)P_{i,n}\right]\right |
\!\!\!\! &\!\!\le\!\!&\!\!\!\!
 \E\left[\left|g(\Xh_T^{t_i,x})-g_\eps(\Xh_T^{t_i,x})P_{i,n}\right|\right]
 \nonumber
 \\
 \!\!\!\!&\!\!\!\!\!\!\!\!&\!\!\!\!
 +\left |\E\left[\left(g_\eps(\Xh_T^{t_i,x})-g_\eps(x)\right)P_{i,n}\right]\right |+|g_\eps-g|_\infty
 \nonumber
 \\
\!\!\!\! &\!\!\!\!\!\!\le\!\!\!\!&\!\!\!\!
 C\eps+\left |\E\left[\left(g_\eps(\Xh_T^{t_i,x})-g_\eps(x)\right)P_{i,n}\right]\right|
 \nonumber
 \\
 \!\!\!\!&\!\!\!\!\!\!\le\!\!\!\!&\!\!\!\!
 C\eps
\!\! +\!\!\left|\E\Big[P_{i,n}\!\!\!\int_{t_i}^T\!\!\!\!\Bigl (\!Dg_\eps\hat b\!+\!\!\frac12\Tr{D^2g_\eps\hat a} \!\Bigr)(s,\Xh_s^{t_i,x})ds\Big]
 \right| 
 \nonumber
 \\
 &&
 +\left|\E\big[P_{i,n}\int_{t_i}^T Dg_\eps(\Xh_s^{t_i,x})\hat\sigma(s) dW_s\big]\right|,
 \label{estim12}
 \ee
where we denoted $\hat b(s)=b(t_j,\Xh_{t_j}^{t_i,x})$ and $\hat\sigma(s)=\sigma(t_j,\Xh_{t_j}^{t_i,x})$ for $t_j\le s<t_{j+1}$ and $\hat a=\hat\sigma^T\hat\sigma$. We next estimate each term separately.  

{\bf 2.a.} First, since $\{P_{i,k},i\le k\le n\}$ is a martingale: 
 \begin{align*}
 \Big|\E\big[P_{i,n}\int_{t_i}^T 
 & Dg_\eps(\Xh_s^{t_i,x})\hat\sigma(s) dW_s\big]\Big|
                 \;=\; \Big|\sum_{j=i}^{n-1}\E\big[P_{i,n}\int_{t_j}^{t_{j+1}} Dg_\eps(\Xh_s^{t_i,x})\hat\sigma(s) dW_s\big]\Big|
 \\
 &\le\;
 \sum_{j=i}^{n-1}
 \Big|\E\big[P_{i,j+1}\int_{t_j}^{t_{j+1}} Dg_\eps(\Xh_s^{t_i,x})\hat\sigma(s) dW_s\big]\Big|
 \\
 &=\;
 \sum_{j=i}^{n-1}
 \Big|\E\big[P_{i,j}\hat\sigma(t_j)\E_{t_j}\big[P_{j,j+1}\int_{t_j}^{t_{j+1}}Dg_\eps(\Xh_s^{t_i,x}) dW_s\big]\big]\Big|.
 \end{align*}
Notice that 
\b*
\E_{t_j}\left[P_{j,j+1}\int_{t_j}^{t_{j+1}}\!\!\!\!\!\!Dg_\eps(\Xh_s^{t_i,x}) dW_s\right]
\!\! &\!\!=\!\!&\!\!
 \E_{t_j}\left[(W_{t_{j+1}}-W_{t_j})^2\int_{t_j}^{t_{j+1}}\!\!\!\!\!\!Dg_\eps(\Xh_s^{t_i,x}) dW_s\right] \\
\!\!&\!\!=\!\!&\!\! \E_{t_j}\left[\int_{t_j}^{t_{j+1}}2W_sDg_\eps(\Xh_s^{t_i,x}) ds\right].
\e*
 Using Lemma \ref{lemintpart} and \reff{geps-g}, this provides: 
 \begin{align}
 \Big|\E\big[P_{i,n}\int_{t_i}^T
 & Dg_\eps(\Xh_s^{t_i,x})\hat\sigma(s) dW_s\big]\Big|
 \\
 &\le\;
2 \sum_{j=i}^{n-1}
 \Big|\E\big[P_{i,j+1}\hat\sigma(t_j)^2\frac{\alpha_j}{h}
             \E_{t_j}\big[\int_{t_j}^{t_{j+1}}sD^2g_\eps(\Xh_s^{t_i,x})ds
                     \big]
         \big]\Big|,
 \nonumber
 \\
 &\le\;
 C\eps^{-1}\sum_{j=i}^{n-1}h
 \;\le\; 
 C'(T-t_i)\eps^{-1}.
 \label{estim1}
 \end{align}

{\bf 2.c.} By \reff{geps-g} and the boundedness of $b$ and $\sigma$, we also estimate that:
\be\label{estim2}
\left |Dg_\eps(\Xh_s^{t_i,x})\hat b(s,\Xh_s^{t_i,x})\!\!+\!\!\frac12\Tr{D^2g_\eps(\Xh_s^{t_i,x})\hat a(s,\Xh_s^{t_i,x})}\right | 
\!\!\!&\!\!\le\!\!&\!\!\!
C+C\eps^{-1}.
\ee

{\bf 2.b.} Plugging \reff{estim1} and \reff{estim2} into \reff{estim12}, we obtain:
 \b*
 \left |\E\left[\left(g_\eps(\Xh_T^{t_i,x})-g_\eps(x)\right)P_{i,n}\right]\right |
 &\le&
 C(T-t_i)+C(T-t_i)\eps^{-1},
 \e*
which by \reff{inequality} provides:
 \b*
 |v^h(t_i,x)-g(x)|
 &\le&
 C\eps+C(T-t_i)\eps^{-1}+C(T-t_i).
 \e*
The required result follows from the choice $\eps=\sqrt{T-t_i}$. 
 \ep

\begin{Corollary}\label{prophalfhold}
The function $v^h$ is 1/2-H\"older continuous on $t$ uniformly on $h$.
\end{Corollary}

\proof
The proof of $\frac12$-H\"older continuity with respect to $t$ could be easily provided by replacing $g$ and $v^h(t_k,\cdot)$ in the assertion of Lemma respectively by $v^h(t,\cdot)$ and $v^h(t',\cdot)$ and consider the scheme from $0$ to time $t'$ with time step equal to $h$. Therefore, we can write;
\b*
|v^h(t,x)-v^h(t',x)|
&\le&
C(t'-t)^{\frac12},
\e*
where $C$ could be chosen independent of $t'$ for $t'\le T$. 
\ep

\subsection{Derivation of the rate of convergence}
\label{proofrate}

The proof of Theorem \ref{thmrateconv} is based on Barles and Jakobsen \cite{barlesjakobsen}, which uses switching systems approximation and the Krylov method of shaking coefficients \cite{krylov}. 

\subsubsection{Comparison result for the scheme}

Because $F$ does not satisfy the standard monotonicity condition \reff{standardmonotonicity} of Barles and Souganidis \cite{barlessouganidis}, we need to introduce the nonlinearity $\overline{F}$ of Remark \ref{remFbar} so that $\overline{F}$ satisfies \reff{strict}. Let $u^h$ be the familiy of functions defined by
 \be\label{defuh}
 u^h(T,.)=g
 ~\mbox{and}~
 u^h(t_i,x)=\overline{\TT}_h [u^h](t_i,x),
 \ee
where for a function $\psi$ from $[0,T]\times\R^d$ to $\R$ with exponential growth:
 \b*
 \overline{\TT}_h[\psi](t,x):=\E\left[\psi(t+h,\hat X_h^{t,x})\right]+h\overline{F}\left(\cdot,\Dc_h\psi\right)(t,x),
 \e*
and set
 \be\label{defbarv}
 \vo^h(t_i,x) \;:=\; e^{-\theta(T-t_i)}u^h(t_i,x),
 &&
 i=0,\ldots,n.
 \ee
The following result shows that the difference $v^h-\vo^h$ is of higher order, and thus reduces the error estimate problem to the analysis of the difference $\vo^h-v$.

\begin{Lemma}\label{lemvbarv}
Under Assumption $\FF$, we have
 \b*
 \limsup_{h\searrow 0} h^{-1}| (v^h-\vo^h)(t,.)|_\infty
 &<&
 \infty.
 \e*
\end{Lemma}

\proof
By definition of $\overline{F}$, we directy calculate that:
 \b*
 \vo^h(t,x)
 &=&
 e^{-\theta h}(1+h\theta)\E[\vo^h(t+h,\hat X^{t,x}_h)]
 +h F\left(t+h,x,\Dc_h\vo^h(t,x)\right).
 \e*
Since $1+h\theta=e^{\theta h}+O(h^2)$, this shows that $\vo^h(t,x)=\TT_h[\vo^h](t,x)+O(h^2)$. By lemma \ref{lemstab}, we conclude that:
\b*
|(\vo^h-v^h)(t,\cdot)|_\infty
&\le&
(1+Ch)|(\vo^h-v^h)(t+h,\cdot)|_\infty+O(h^2),
\e*
which shows by the Gronwall inequality that $|(\vo^h-v^h)(t,\cdot)|_\infty
\le O(h)$ for all $t\le T-h$.
 \ep 
 
\vspace{5mm}

By Remark \ref{remFbar}, the operator $\overline{\TT}_h$ satisfies the standard monotonicity condition \reff{standardmonotonicity}:
 \be\label{standardmonotonicitybar}
 \varphi \;\le\; \psi
 &\Longrightarrow& 
 \overline{\TT}_h[\varphi] \;\le\; \overline{\TT}_h[\psi].
 \ee
The key-ingredient for the derivation of the error estimate is the following comparison result for the scheme.

\begin{Proposition} \label{propmaxpri}
Let Assumption $\FF$ holds true, and set $\beta:=|F_r|_\infty$. Consider two arbitrary bounded functions $\varphi$ and $\psi$ satisfying:
 \be\label{hyplemmaxpri}
 h^{-1}\left(\varphi-\overline{\TT}_h[\varphi]\right) \;\le\; g_1
 &\mbox{and}&
 h^{-1}\left(\psi-\overline{\TT}_h[\psi]\right) \;\ge\; g_2
 \ee
for some bounded functions $g_1$ and $g_2$. Then, for every $i=0,\cdots,n$:
 \be\label{conlemmaxpri}
 (\varphi-\psi)(t_i,x)
 &\le&
 e^{\beta(T-t_i)}|(\varphi-\psi)^+(T,\cdot)|_\infty
 +(T-h)e^{\beta(T-t_i)}|(g_1-g_2)^+|_\infty.
 \ee
\end{Proposition}

To prove this comparison result, we need the following strengthening of the  monotonicity condition:

\begin{Lemma} \label{lemratemon}
Let Assumption $\FF$ hold true and let $\beta:=|F_r|_\infty$. Then, for every $a,b\in\R_+$, and every bounded functions $\varphi\le \psi$, the function $\delta(t):=e^{\beta(T-t)}(a+b(T-t))$ satisfies:
 \b*
 \overline{\TT}_h[\varphi+\delta](t,x)
 &\le& 
 \overline{\TT}_h[\psi](t,x)+\delta(t)-hb,
 ~~
 t\le T-h,~x\in\R^d.
 \e*
\end{Lemma}
\proof
Because $\delta$ does not depend on $x$, we have $\Dc_h[\varphi+\delta]=\Dc_h\varphi+\delta(t+h)e_1$, where $e_1:=(1,0,0)$. Then, it follows from 
the regularity of $\overline{F}$ that there exist some $\xi$ such that:
 \begin{align*}
 \overline{F}\big(t+h,x,\Dc_h[\varphi+\delta](t,x)\big)
 =&
 \overline{F}\big(t+h,x,\Dc_h\varphi(t,x)\big)\\
& +\delta(t+h)\overline{F}_r\big(t+h,x,\xi e_1+\Dc_h\varphi(t,x)\big),
 \end{align*}
and
 \b*
 \overline{\TT}_h[\varphi+\delta](t,x)
\!\! &\!\!=\!\!&\!\!
 \delta(t+h)+\E[\varphi(t+h,\hat X_h^{t,x})]
 +h\overline{F}\big(t+h,x,\Dc_h\varphi(t,x)\big)
 \\ 
\!\!&\!\!\!\! &\!\!
 + h\delta(t+h)\overline{F}_r\big(t+h,x,\xi e_1+\Dc_h\varphi(t,x)\big)
 \\
\!\! &\!\!=\!\!&\!\! 
 \overline{\TT}_h[\varphi](t,x)
 +\delta(t+h)\left\{1+h\overline{F}_r\big(t+h,x,\xi e_1+\Dc_h\varphi(t,x)\big)\right\}
 \\
\!\! &\!\!\le\!\!&\!\!
 \overline{\TT}_h[\varphi](t,x)+\left(1+\beta h\right)\delta(t+h).
 \e*
Since $\overline{\TT}_h$ satisfies the standard monotonicity condition \reff{standardmonotonicitybar}, this provides:
 \b*
 \overline{\TT}_h[\varphi\!+\!\delta](t,x)
\!\! \le\!\!
 \overline{\TT}_h[\psi](t,x)\!+\!\delta(t)\!+\!\zeta(t),
 &\!\!\mbox{where}\!\!&
 \zeta(t) 
 \!\!:=\!\!
 \left(1\!+\!\beta h\right)\delta(t\!+\!h)\!-\!\delta(t).
 \e*
It remains to prove that $\zeta(t)$ $\le-hb$. From the smoothness of $\delta$, we have $\delta(t+h)-\delta(t)=h\delta'(\bar t)$ for some $\bar t\in[t,t+h)$. Then, since $\delta$ is decreasing in $t$, we see that
 \b*
 h^{-1}\zeta(t)
 &=&
 \delta'(\bar t)+\beta \delta(t+h)
 \;\le\;
 \delta'(\bar t)+\beta \delta(\bar t)
 \;\le\;
 -b e^{\beta(T-\bar t)},
 \e*
and the required estimate follows from the restriction $b\ge 0$.
 \ep

\vspace{5mm}

\no {\bf Proof of Proposition \ref{propmaxpri}.}\quad We may refer directly to the similar result of \cite{barlesjakobsen}. However in our context, we give the following simpler proof. Observe that we may assume without loss of generality that 
 \be\label{wlog}
 \varphi(T,\cdot)\le \psi(T,\cdot)
 &\mbox{and}&
 g_1\le g_2.
 \ee 
Indeed, one can otherwise consider the function
 \b*
 \bar\psi
\! :=\!
 \psi\!+\!e^{\beta(T-t)}\!\left(a\!+\!b(T\!-\!t)\right)
 \!\!&\mbox{where}&\!\!
 a\!=\!|(\varphi-\psi)^+(T,\cdot)|_\infty,
 ~
 b\!=\!|(g_1\!-\!g_2)^+|_\infty,
 \e*
and $\beta$ is the parameter defined in the previous Lemma \ref{lemratemon}, so that $\bar\psi(T,\cdot)\ge\varphi(T,\cdot)$ and, by Lemma \reff{lemratemon}, $\bar\psi(t,x)-\overline{\TT}_h[\bar\psi](t,x)\ge h (g_1\vee g_2)$. 
Hence \reff{wlog} holds true for $\varphi$ and $\bar\psi$.

We now prove the required result by induction. First $\varphi(T,\cdot)\le\psi(T,\cdot)$ by \reff{wlog}. We next assume that $\varphi(t+h,\cdot)\le\psi(t+h,\cdot)$ for some $t+h\le T$. Since $\overline{\TT}_h$ satisfies the standard monotonicity condition \reff{standardmonotonicitybar}, it follows from \reff{wlog} that
 \b*
 \overline{\TT}_h[\varphi](t,x)
 &\le&
 \overline{\TT}_h[\psi](t,x).
 \e*
On the other hand, under \reff{wlog}, the hypothesis of the lemma implies:
 \b*
 \varphi(t,x)-\overline{\TT}_h[\varphi](t,x)
 &\le&
 \psi(t,x)-\overline{\TT}_h[\psi](t,x).
 \e*
Then $ \varphi(t,\cdot)\le\psi(t,\cdot)$.
 \ep

\subsubsection{Proof of Theorem \ref{thmrateconv} (i)}

Under the conditions of Assumption HJB on the coefficients, we may build a bounded subsolution $\underline{v}^\eps$ of the nonlinear PDE, by the method of shaking the coefficients, which is Lipschitz in $x$, $1/2-$H\"older continuous in $t$, and approximates uniformly the solution $v$:
 $$
 v-\eps \le v^\eps \le v.
 $$
Let $\rho(t,x)$ be a $C^\infty$ positive function supported in $\{(t,x):t\in[0,1],|x|\le 1\}$ with unit mass, and define
 \be\label{defoverueps}
 \underline{w}^\eps(t,x)
 :=v^\eps*\rho^\eps
 &\mbox{where}&
 \rho^\eps(t,x)
 :=\frac{1}{\eps^{d+2}}\rho\left(\frac{t}{\eps^2},\frac{x}{\eps}\right)
 \ee 
so that, from the convexity of the operator $F$,
 \be\label{vunderlineweps}
 \underline{w}^\eps
 &\mbox{is a subsolution of \reff{equation},}&
 |\underline{w}^\eps-v|\le 2\eps.
 \ee
Moreover, since $v^\eps$ is Lipschitz in $x$, and $1/2-$H\"older continuous in $t$, 
 \be\label{regunderweps}
 \underline{w}^\eps~
 \mbox{is}~C^\infty,
 &\mbox{and}&
 \left|\partial_t^{\beta_0}D^\beta\underline{w}^\eps\right| 
 \le C\eps^{1-2\beta_0-|\beta|_1}\\
\nonumber &\mbox{for any}&(\beta_0,\beta)\in\N\times\N^d\setminus\{0\},
 \ee
where $|\beta|_1:=\sum_{i=1}^d\beta_i$, and $C>0$ is some constant.
As a consequence of the consistency result of Lemma \ref{lemcons} above, we know that
 \b*
 \Rc_h[\underline{w}^\eps](t,x)
\!\! &\!\!\!:=\!\!\!&\!\!
 \frac{\underline{w}^\eps(t,x)\!-\!\TT_h[\underline{w}^\eps](t,x)}{h}
 \!+\!
 \mathcal{L}^X\underline{w}^\eps(t,x)\!+\!F(\cdot,\underline{w}^\eps,D\underline{w}^\eps,D^2\underline{w}^\eps)(t,x)
 \e*
converges to $0$ as $h\to 0$. The next key-ingredient is to estimate the rate of convergence of $\Rc_h[\underline{w}^\eps]$ to zero:

\begin{Lemma}\label{lemratcons}
For a family $\{\varphi_\eps\}_{0<\eps<1}$ of smooth functions satisfying \reff{regunderweps}, we have: 
 \b*
 \left|\Rc_h[\varphi_\eps]\right|_\infty
 \;\le\;
 R(h,\eps):=C\;h\eps^{-3}
 &\mbox{for some constant}& C>0.
 \e*
\end{Lemma}

The proof of this result is reported at the end of this section. From the previous estimate together with the subsolution property of $w^\eps$, we see that $\underline{w}^\eps \le \overline{\TT}_h[\underline{w}^\eps]+ Ch^2\eps^{-3}$. Then, it follows from Proposition \ref{propmaxpri} that
 \be\label{underlineww}
 \underline{w}^\eps-\vo^h\le C|(\underline{w}^\eps-\vo^h)(T,.)|_\infty+Ch\eps^{-3}\le C(\eps+h\eps^{-3}).
 \ee
We now use \reff{vunderlineweps} and \reff{underlineww} to conclude that
 \b*
 v-\vo^h
 \le
 v-\underline{w}^\eps+\underline{w}^\eps-\vo^h
 \le C(\eps + h\eps^{-3}).
 \e*
Minimizing the right hand-side estimate over the choice of $\eps>0$, this implies the upper bound on the error $v-v^h$:
 \be\label{upperbound}
 v-\vo^h 
 &\le& 
 Ch^{1/4}.
 \ee

\subsubsection{Proof of Theorem \ref{thmrateconv} (ii)}

The results of the previous section, together with the reinforced assumption HJB+, allow to apply the switching system method of Barles and Jakobsen \cite{barlesjakobsen} which provides the lower bound on the error:
 \b*
 v-\vo^h &\ge& -\inf_{\eps>0}\{C\eps^{1/3}+R(h,\eps)\}
 \;=\;
 -C'h^{1/10},
 \e*
for some constants $C,C'>0$. The required rate of convergence follows again from Lemma \ref{lemvbarv} which states that the difference $v^h-\overline{v}^h$ is dominated by the above rate of convergence.

\vspace{5mm}

\no {\bf Proof of Lemma \ref{lemratcons}}\quad Notice that the evolution of the Euler approximation $\hat X^{t,x}_h$ between $t$ and $t+h$ is driven by a constant drift $\mu(t,x)$ and a constant diffusion $\sigma(t,x)$. Since $D\varphi_\eps$ is bounded, it follows from It\^o's formula that:
 $$
 \frac{1}{h}\!\!\left[\E\varphi_\eps(t\!+\!h,\hat X_h^x)\!-\!\varphi_\eps(t,x)\right]\!-\!\Lc^X\!\varphi_\eps(t,x)
\!\! =\!\!
 \frac{1}{h}\E\!\!\int_t^{t+h}\!\!\!\!\!\!\!\! \bigl(\Lc^{\hat X^{t,x}}\!\!\!\!\varphi_\eps(u,\hat X_u^x)-\Lc^X\!\varphi_\eps(t,x)\!\bigr)\!du,
 $$ 
where $\Lc^{\hat X^{t,x}}$ is the Dynkin operator associated to the Euler scheme:
 \b*
 \Lc^{\hat X^{t,x}}\varphi(t',x')
 &=&
 \partial_t\varphi(t',x')+\mu(t,x)D\varphi(t',x')+\frac12{\rm Tr}\left[a(t,x)D^2\varphi(t',x')\right].
 \e*
Applying again It\^o's formula, and using the fact that $\Lc^{\hat X^{t,x}}D\varphi_\eps$ is bounded, leads to
 \b*
 \frac{1}{h}\!\left[\E\varphi_\eps(t\!+\!h,\hat X_h^x)\!-\!\varphi_\eps(t,x)\right]\!-\!\Lc^X\!\varphi_\eps(t,x)
\!\! &\!\!\!=\!\!\!&\!\!
 \frac{1}{h}\E\!\!\int_t^{t+h}\!\!\!\! \int_t^u\!\!\!\Lc^{\hat X^{t,x}}\!\!\Lc^{\hat X^{t,x}}\!\!\!\!\!\varphi_\eps(s,\hat X_s^x)dsdu.
 \e*
Using the boundedness of the coefficients $\mu$ and $\sigma$, it follows from \reff{regunderweps} that for $\eps\in(0,1)$:
 \b*
 \left|\frac{\E\varphi_\eps(t+h,\hat X_h^x)-\varphi_\eps(t,x)}{h}-\Lc^X\varphi_\eps(t,x)\right| 
 &\le&
 R_0(h,\eps):=C\;h\eps^{-3}.
 \e*
\no\underline{\it Step 2}\quad This implies that
 \be
 \left|\Rc_h[\varphi_\eps](t,x)\right|
 \!\!&\!\!\le\!\!&\!\!
 \left|\frac{\E\varphi_\eps(t+h,\hat X_h^{t,x})-\varphi_\eps(t,x)}{h}-\Lc^X\varphi_\eps(t,x)\right|
 \nonumber\\
 \!\!&\!\!\!\!&\!\!
 +
 \left|F(x,\varphi_\eps(t,x),D\varphi_\eps(t,x),D^2\varphi_\eps(t,x))-F\left(\cdot,\Dc_h[\varphi_\eps](t,x)\right)\right|
 \nonumber\\
\!\! &\!\!\le\!\!&\!\!
 R_0(h,\eps)+ C\sum_{k=0}^2 \left|\E D^k\varphi_\eps(t+h,\hat X_h^{t,x})-D^k\varphi_\eps(t,x)\right|
 \label{step1}
 \ee
by the Lipschitz continuity of the nonlinearity $F$.

By a similar calculation as in Step 1, we see that: 
 \b*
 |\E D^i\varphi_\eps(t+h,\hat X_h^{t,x})-D\varphi_\eps(t,x)|
 &\le&
 Ch\eps^{-1-i},~i=0,1,2,
 \e*
which, together with \reff{step1}, provides the required result.
 \ep

\vspace{5mm}
\subsection{The rate of convergence in the linear case}
\label{subsectlinear}

In this subsection, we specialize the discussion to the linear one-dimensional case
 \be\label{Flinear}
 F(\gamma)
 &=&
 c\gamma,
 \ee
for some $c>0$. The multi-dimensional case $d>1$ can be handled similarly. Assuming that $g$ is bounded, the linear PDE \reff{equation}-\reff{terminal} has a unique bounded solution
 \be\label{vExp}
 v(t,x)
 \;=\;
 \E\left[ g\left(x+\sqrt{1+2c}\;W_{T-t}\right) \right]
 &\mbox{for}&
 (t,x)\in[0,T]\times\R^d.
 \ee
We also observe that this solution $v$ is $C^\infty\left([0,T)\times\R\right)$ with
 \be\label{Dkv}
 D^{k}v(t,x)
 \;=\;
 \E\left[ g^{(k)}\left(x+\sqrt{1+2c}\;W_{T-t}\right) \right],
 &t<T,& x\in\R.
 \ee
This shows in particular that $v$ has bounded derivatives of any order, whenever the terminal data $g$ is $C^\infty$ and has bounded derivatives of any order. 

Of course, one can use the classical Monte Carlo estimate to produce an approximation of the function $v$ of \reff{vExp}. The objective of this section is to analyze the error of the numerical scheme outlined in the previous sections. Namely:
 \begin{equation}
 v^h(T,\cdot)\!=\!g,
 ~
 v^h(t_{i-1},x)
 \!=\!
 \E\!\left[v^h(t_i,x+W_h)\right]
 \!+ ch\E\!\left[v^h(t_i,x+W_h)H^h_2\right]\!\!,
 ~i\le n.
 \end{equation}
Here, $\sigma=1$ and $\mu=0$ are used to write the above scheme.

\begin{Proposition}\label{proplinear}
Consider the linear $F$ of \reff{Flinear}, and assume that $D^{(2k+1)}v$ is bounded for every $k\ge 0$. Then
 \b*
 \limsup_{h\to 0} h^{-1/2}|v^h-v|_\infty &<& \infty.
 \e*
\end{Proposition}

\proof Since $v$ has bounded first derivative with respect to $x$, it follows from It\^o's formula that:
 \b*
 v(t,x)
 &=&
 \E\left[v(t+h,x+W_h)\right] + c\E\left[\int_0^h \lap v(t+s,v+W_s)ds\right],
 \e*
Then, in view of Lemma \ref{lemintpart}, the error $u:=v-v^h$ satisfies $u(t_n,X_{t_n})=0$ and for $i\le n-1$:
 \be
 u\left(t_i,X_{t_i}\right)
 &=&
 \E_i\left[u\left(t_{i+1},X_{t_{i+1}}\right)\right]
 +ch\;\E_i\left[\lap u\left(t_{i+1},X_{t_{i+1}}\right)\right]\nonumber\\
 &&
 +c\E_i\int_0^h \left[\lap v\left(ih+s,X_{ih+s}\right)-\lap v\left((i+1)h,X_{(i+1)h}\right)\right]ds,
 \label{xstart}
 \ee
where $\E_i:=\E[\cdot|\mathcal{F}_{t_i}]$ is the expectation operator conditional on $\Fc_{t_i}$.

\no\underline{\it Step 1}\quad 

Set
 \b*
 a^k_i \!\!:=\!\! \E\!\left[\lap^k u\left(t_i,X_{t_i}\right)\right]\!,~
 b^k_i \!\!:=\! \E\!\int_0^h \!\left[\lap^k v\left(t_{i-1}\!+\!s,X_{t_{i-1}\!+\!s}\right)\!-\!\lap^k v\left(t_i,X_{t_i}\right)\right]ds,
 \e*
and we introduce the matrices
 \b*
 A:=\left(\begin{array}{ccccc} 1 & -1    & 0      & \cdots &    0 \\
                               0 & 1     & -1     & \cdots &    0 \\
                          \vdots &\ddots & \ddots & \ddots &\vdots \\
                          \vdots & \ddots & \ddots &   1    & -1\\
                               0 & \cdots & \cdots &   0    & 1
          \end{array}
    \right),
 &&
 B:=\left(\begin{array}{ccccc} 0 & 1    & 0    &\ldots& 0 \\
                          \vdots &\ddots&\ddots&\ddots& \vdots\\
                          \vdots &\ddots&\ddots&\ddots& 0\\
                          \vdots &\ddots&\ddots&\ddots& 1\\
                               0 &\cdots&\cdots&\cdots& 0\\
          \end{array}
    \right),
 \e*
and we observe that \reff{xstart} implies that the vectors $a^k:=(a^k_1,\ldots,a^k_n)^{\rm T}$ and $b^k:=(b^k_1,\ldots,b^k_n)^{\rm T}$ satisfy $A a^k = ch B a^{k+1} + c B b^k$ for all $k\ge 0$, and therefore:
 \begin{equation}\label{recurrence}
 a^k \!=\! ch A^{-1}B a^{k+1} \!+\! c A^{-1}B b^k \!,
~\mbox{where} ~
 A^{-1}\!=\!\left(\begin{array}{cccc} 1 & 1    & \cdots &    1 \\
                                  0 & 1    & \cdots &    1 \\
                             \vdots &\ddots& \ddots &\vdots \\
                                  0 & \cdots & 0    & 1
              \end{array}
        \right).
 \end{equation}
By direct calculation, we see that the powers $(A^{-1}B)^{k}$ are given by:
 \b*
 (A^{-1}B)^{k}_{i,j}
 =
 \1_{\{j\ge i+k\}}{{j-i-1}\choose{k-1}}
 &\mbox{for all}&
 k\ge 1~\mbox{and}~i,j=1,\ldots,n.
 \e*
In particular, because $a^k_n=0$, $(A^{-1}B)^{n-1}a^k=0$. Iterating \reff{recurrence}, this provides:
 \b*
 a^0 &=& ch(A^{-1}B)a^1+c(A^{-1}B)b^0 \;=\;\ldots\;=\;\sum_{k=0}^{n-2}c^{k+1} h^k(A^{-1}B)^{k+1}b^k,
 \e*
and therefore:
 \be\label{errordecomplinear}
 u(0,x)
 &=&
 a^0_1
 \;=\;
 c\sum_{k=0}^{n-2} (ch)^k (A^{-1}B)^{k+1}_{1,j} b^k. 
 \ee
Because of 
 \b*
 (A^{-1}B)^{k}_{1,j}
 =
 \1_{\{j\ge 1+k\}}{{j-2}\choose{k-1}}
 &\mbox{for all}&
 k\ge 1~\mbox{and}~j=1,\ldots,n \; ,
 \e*
 we can write \reff{errordecomplinear}:
 \b*
 u(0,x)
 &
 \;=\;
 &
 c\sum_{k=0}^{n-2} (ch)^k\sum_{j=k+2}^n{{j-2}\choose{k}} b^{k-1}_j. 
 \e*
 By changing the order of the summations in the above we conclude that:
 \be
 \label{errordecomplinear1}
 u(0,x)
 &
 \;=\;
 &
 c\sum_{j=2}^{n} \sum_{k=0}^{j-2} (ch)^k{{j-2}\choose{k}} b^{k-1}_j. 
 \ee 
 \no\underline{\it Step 2}\quad From our assumption that $D^{2k+1}v$ is $\L^\infty-$bounded for every $k\ge 0$, it follows that
 \b*
 |b^k_j|
 \;\le\;
 \E\left[\int_{t_{i-1}}^{t_i} \left|\lap^k v(s,X_s)-\lap^k v(t_j,X_{t_j})\right|ds\right]
 &\le& C h^{3/2}
 \e*
for some constant $C$. We then deduce from \reff{errordecomplinear1} that:
 \b*
 |u(0,x)|
 &\le&  c Ch^{3/2}\;\sum_{j=2}^{n} \sum_{k=0}^{j-2} (ch)^k{{j-2}\choose{k}}.
 \e*
 So,
 \b*
 |u(0,x)|
 &\le&
 c Ch^{3/2}\;\sum_{j=2}^n (1+ch)^{j-2}
  \;=\;
 c Ch^{3/2}\;\frac{(1+ch)^{n-1}-1}{ch}
 \;\le\;
 C\; \sqrt{h}.
 \e* 

\section{Probabilistic Numerical Scheme}
\setcounter{equation}{0}
\label{Sec:sectprobnum}

In order to implement the backward scheme \reff{scheme}, we still need to discuss the numerical computation of the conditional expectations involved in the definition of the operators $\TT_h$ in \reff{TT}. In view of the Markov feature of the process $X$, these conditional expectations reduce to simple regressions. Motivated by the problem of American options in financial mathematics, various methods have been introduced in the literature for the numerical approximation of these regressions. We refer to \cite{bt} and \cite{glw} for a detailed discussion.

The chief object of this section is to investigate the asymptotic properties of our suggested numerical method when the expectation operator $\E$ in \reff{scheme} is replaced by some estimator $\hat \E^N$ corresponding to a sample size $N$:
\be
 \label{operator}
 \tilde \TT_h^N[\psi](t,x)
 &
 :=
 &
 \hat\E^N\left[\psi(t+h,\hat X_h^x)\right] + hF\left(\cdot,\hat\Dc_h\psi\right)(t,x),\\
 \label{trunc}
 \hat \TT_h^N[\psi](t,x)
 &
 :=
 &
 -K_h[\psi]\vee\tilde \TT_h^N[\psi](t,x)\wedge K_h[\psi]
 \ee
where 
\b*
\hat \Dc_h\psi(t,x)
 :=
\hat\E^N\left[\psi(t+h,\hat X^{t,x}_h)H_h(t,x)\right], 
 ~K_h[\psi]
 :=
\Vert \psi\Vert_\infty(1+C_1h)+C_2h,
\e*
where
\b*
C_1=\frac{1}{4}\vert F_p^\text{T} F_\gamma^- F_p\vert_\infty+\vert F_r\vert_\infty
&\mbox{and}&
C_2=\vert F(t,x,0,0,0)\vert_\infty.
\e*
The above bounds are needed for technical reasons which were already observed in \cite{bt}. 

With these notations, the implementable numerical scheme is:
 \be
 \label{schememc}
 \hat v^h_{N}(t,x,\omega)=\hat\TT_h^{N}[\hat v^h_{N}](t,x,\omega),
 \ee
where $\hat\TT_h^{N}$ is defined in \reff{operator}-\reff{trunc}, and the presence of $\omega$ throughout this section emphasizes the dependence of our estimator on the underlying sample.

Let $\Rc_b$ be the family of random variables $R$ of the form $\psi(W_h)H_i(W_h)$ where $\psi$ is a function with $|\psi|_\infty\le b$ and $H_i$'s are the Hermite polynomials:
\b*
H_0(x)=1,\;H_1(x)=x\;\;{\rm and}\;\;H_2(x)=x^Tx-h\;\;\forall x\in R^d.
\e* 

\no {\bf Assumption E}\quad 
{\it
There exist constants $C_b,\lambda,\nu>0$ such that $$\left\Vert\hat\E^N[R]-\E[R]\right\Vert_p \le  C_bh^{-\lb}N^{-\nu}$$ for every
$R\in\Rc_b$, for some $p\ge 1$.}

\vspace{5mm}

\begin{Example}{\rm
Consider the regression approximation based on the Malliavin integration by parts as introduced in Lions and Reigner \cite{lionsr}, Bouchard, Ekeland and Touzi \cite{bet}, and analyzed in the context of the simulation of backward stochastic differential equations by \cite{bt} and \cite{cmt}. Then Assumption {\bf E} is satisfied for every $p>1$ with the constants $\lb=\frac{d}{4p}$ and $\nu=\frac{1}{2p}$, see \cite{bt}.}
\end{Example}

Our next main result establishes conditions on the sample size $N$ and the time step $h$ which guarantee the convergence of $\hat v^h_{N}$ towards $v$. 

\begin{Theorem}\label{thmconvmc}
Let Assumptions {\rm\bf E} and {\rm\bf F} hold true, and assume that the fully nonlinear PDE \reff{equation} has comparison with growth $q$. Suppose in addition that
\be
\label{samplesize1}
\lim_{h\to 0} h^{\lb+2}N_h^\nu
&=&
\infty.
\ee
Assume that the final condition $g$ is bounded Lipschitz, and the coefficients $\mu$ and $\sigma$ are bounded. Then, for almost every $\omega$: 
\b*
 \hat v^h_{N_h}(\cdot,\omega)\longrightarrow v
 &&
 \mbox{locally uniformly,}
\e*
where v is the unique viscosity solution of \reff{equation}.
\end{Theorem}

\proof 
We adapt the argument of \cite{barlessouganidis} to the present stochastic context.
By Remark \ref{remFbar} and Lemma \ref{lemvbarv}, we may assume without loss of generality that the strict monotonicity \reff{strict} holds.

By \reff{trunc}, we see that $\hat v^h$ is uniformly bounded. So, we can define:
\be
\label{liminf}
\hat v_*(t,x):=\mathop{\liminf}\limits_{\tiny{\begin{array}{c}
             (t',x')\to(t,x)
             \\
             h\to 0
             \end{array}}}
\hat v^h(t',x')
 &\mbox{and}&
 \hat v^*(t,x):=\mathop{\limsup}\limits_{\tiny{\begin{array}{c}
             (t',x')\to(t,x)
             \\
             h\to 0
             \end{array}}}
\hat v^h(t',x').
\ee
Our objective is to prove that $\hat v_*$ and $\hat v^*$ are respectively viscosity superpersolution and subsolution of \reff{equation}. By the comparison assumption, we shall then conclude that they are both equal to the unique viscosity solution of the problem whose existence is given by Theorem \ref{thmconv}. In particular, they are both deterministic functions. 

We shall only report the proof of the supersolution property, the subsolution property follows from the same type of argument.

In order to prove that $\hat v_*$ is a supersolution of \reff{equation}, we consider $(t_0,x_0)\in[0,T)\times\R^n$ together with a test function $\varphi\in C^2\left([0,T)\times\R^n\right)$, so that 
\b*
0=\mathop{\min}\{ \hat v_*-\varphi\}=(\hat v_*-\varphi)(t_0,x_0).
\e*
By classical manipulations, we can find a sequence $(t_n,x_n,h_n)\to(t_0,x_0,0)$ so that $\hat v^{h_n}(t_n,x_n)\to \hat v_*(t_0,x_0)$ and
\b*
(\hat v^{h_n}-\varphi)(t_n,x_n)=\mathop{\min}\{\hat v^{h_n}-\varphi\}=:C_n\to0.
\e*
Then, $\hat v^{h_n}\ge \varphi+C_n$, and it follows from the monotonicity of the operator $\TT_h$ that:
\b*
\TT_{h_n} [\hat v^{h_n}]\ge \TT_{h_n} [\varphi+C_n].
\e*
By the definition of $\hat v^{h_n}$ in \reff{schememc}, this provides: 
\b*
\hat v^{h_n}(t,x)\ge  \TT_{h_n} [\varphi+C_n](t,x) -(\TT_{h_n}-\hat \TT_{h_n})[\hat v^h_n](t,x),
\e*
where, for ease of notations, the dependence on $N_h$ has been dropped.
Because $\hat v^{h_n}(t_n,x_n)=\varphi(t_n,x_n)+C_n$, the last inequality gives:
\b*
\varphi(t_n,x_n)+C_n-\TT_{h_n} [\varphi+C_n](t_n,x_n)+
h_n R_n\ge 0,
\e*
where $
R_n:=h_n^{-1}
(\TT_{h_n}-\hat \TT_{h_n})[\hat v^{h_n}](t_n,x_n).$
 
We claim that
 \be\label{Rntozero}
 R_n \longrightarrow 0
 &&
 \P-\mbox{a.s. along some subsequence.}
 \ee
Then, after passing to the subsequence, dividing both sides by $h_n$, and sending $n\to\infty$, it follows from Lemma \ref{lemcons} that:
\b*
-\Lc^X \varphi-F\left(\cdot,\varphi,D \varphi,D^2\varphi\right)\ge 0,
\e*
which is the required supersolution property.

It remains to show \reff{Rntozero}. We start by bounding $R_n$ with respect to the error of estimation of conditional expectation.
By Lemma \ref{lemstab}, $\vert\TT_{h_n}[\hat v^{h_n}]\vert_\infty\le K_{h_n}$ and so by \reff{trunc}, we can write:
\be
 \label{hatTTtilde}
 \left|\left(\TT_{h_n}-\hat \TT_{h_n}\right)[\hat v^{h_n}](t_n,x_n)\right|
 &
 \le
 & 
 \left|\left(\TT_{h_n}-\tilde \TT_{h_n}\right)[\hat v^{h_n}](t_n,x_n)\right|.
\ee
By the Lipschitz-continuity of $F$, we have:
\b*
\left|\left(\TT_{h_n}-\hat \TT_{h_n}\right)[\hat v^{h_n}](t_n,x_n)\right|
&
\le
&
C\left(\Ec_0+h_n\Ec_1+h_n\Ec_2\right).
\e*
where:
\b*
\Ec_i
&
=
&
|(\E-\hat\E)[\hat v^{h_n}(t_n+h_n,X_{h_n}^{x_n})H_i^{h_n}(t_n,x_n)]|
\e*
Therefore,
\b*
\left|\left(\TT_{h_n}-\hat \TT_{h_n}\right)[\hat v^{h_n}](t_n,x_n)\right|
\!\!&\!\!
\le
\!\!&\!\!
C\Bigl(\left|(\E-\hat\E)[R_n^0]\right|+\left|(\E-\hat\E)[R_n^1]\right|\\
\!\!&\!\!\!\!&\!\!\hspace*{27mm}+h_n^{-1}\left|(\E-\hat\E)[R_n^2]\right|\Bigr).
\e*
where $R_n^i=\hat v^{h_n}\bigl(t_n+h_n,x_n+\sigma(x)W_h\bigr)H_i(W_h)$, $i=1,2,3$ and $H_i$ is Hermite polynomial of degree $i$.
This leads the following estimate for the error $R_n$:
\be
\label{errorestim}
\left|R_n\right|
\!\!&\!\!
\le
\!\!&\!\!
\frac{C}{h_n}\left(\left|(\E-\hat\E)[R_n^0]\right|+\left|(\E-\hat\E)[R_n^1]\right|+h_n^{-1}\left|(\E-\hat\E)[R_n^2]\right|\right).
\ee
Because $R_n^i\in\Rc_b$ with bound obtained in Lemma \ref{lemstab} by Assumption {\bf E} we have,:
\b*
\Vert R_n\Vert_p
&
\le
& 
Ch_n^{-\lb-2}N_{h_n}^{-\nu},
\e*
so by \reff{samplesize1} we have $\Vert R_n\Vert_p\longrightarrow 0$ which implies \reff{Rntozero}.
\ep

\vspace{5mm}

We finally discuss the choice of the sample size so as to keep the same rate for the error bound.

\begin{Theorem}\label{thmratmc}
Let the nonlinearity $F$ be as in Assumption {\rm\bf HJB}, and consider a regression operator satisfying Assumption {\rm \bf E}. Let the sample size $N_h$ be such that
 \be\label{samplesize2}
 \lim_{h\to 0} h^{\lb+\frac{21}{10}}N_h^\nu
 &>&
 0.
 \ee
Then, for any bounded Lipschitz final condition $g$, we have the following $\L^p-$bounds on the rate of convergence:
 \b*
 \Vert v-\hat v^h\Vert_p
 &\le&
 Ch^{1/10}.
 \e*
\end{Theorem}

\proof By Remark \ref{remFbar} and Lemma \ref{lemvbarv}, we may assume without loss of generality that the strict monotonicity \reff{strict} holds true.

We proceed as in the proof of Theorem \ref{thmrateconv} to see that
 \b*
 v-\hat v^h &\le& v-v^h+v_h-\hat v^h \;=\; \eps+R(h,\eps)+v_h-\hat v^h.
 \e*
Since $\hat v^h$ satisfies \reff{schememc},
 \b*
 h^{-1}\left(\hat v^h-\TT_h[\hat v^h]\right)
 \;\ge\;
 -R_h[\hat v^h]
 &\mbox{where}&
 R_h[\varphi]
 \;:=\;
 \frac{1}{h}\left|\left(\TT_h-\hat\TT_h\right)[\varphi]\right|,
 \e*
where, in the present context, $R_h[\hat v^h]$ is a non-zero stochastic term. 
By Proposition \ref{propmaxpri}, it follows from the last inequality that:

 \b*
 v-\hat v^h\le C\left(\eps + R(h,\eps) + R_h[\hat v^h]\right),
 \e*
where the constant $C>0$ depends only on the Lipschitz coefficient of $F$, $\beta$ in Lemma \ref{lemratemon} and the constant in Lemma \ref{lemratcons}.

Similarly, we follow the line of argument of the proof of Theorem \ref{thmrateconv} to show that a lower bound holds true, and therefore:
 \b*
 |v-\hat v^h|\le C\left(\eps^{1/3} + R(h,\eps) + R_h[\hat v^h]\right),
 \e*
We now use \reff{samplesize2} and proceed as in the last part of the proof of Theorem \ref{thmconvmc} to deduce from \reff{errorestim} and Assumption {\bf F} that
\b*
\Vert R_h[\hat v^h]\Vert_p &\le& Ch^{1/10}.
\e*
With this choice of the sample size $N$, the above error estimate reduces to  
\b*
 \Vert\hat v^h-v\Vert_p
 &\le& 
 C\left(\eps^{1/3} + R(h,\eps) + h^{1/10}\right),
\e*
and the additional term $h^{1/10}$ does not affect the minimization with respect to $\eps$.
\ep

\begin{Example}
{\rm Let us illustrate the convergence results of this section in the context of the Malliavin integration by parts regression method of \cite{lionsr} and \cite{bt} where $\lb=\frac{d}{4p}$ and $\nu=\frac{1}{2p}$ for every $p>1$. So, for the convergence result we need to choose $N_h$ of the order of $h^{-\alpha_0}$ with $\alpha_0>\frac{d}{2}+4p$. For the $L^p$-rate of convergence result, we need to choose $N_h$ of the order of $h^{-\alpha_1}$ with $\alpha_1\ge\frac{d}{2}+\frac{21p}{5}$.
}
\end{Example}
\section{Numerical Results}\label{Sec:sectnulmerics}
\setcounter{equation}{0}

In this section, we provide an application of the Monte Carlo-finite differences scheme suggested in this paper in the context of two different types of problems. We first consider the classical mean curvature flow equation as the simplest front propagation example. We test our backward probabilistic scheme on the example where the initial data is given by a sphere, for which an easy explicit solution is available. A more interesting geometric example in space dimensions 2 is also considered. We next consider the Hamilton-Jacobi-Bellman equation characterizing the classical optimal investment problem in financial mathematics. Here, we again test our scheme in dimension two where an explicit solution is available, and we consider more involved examples in space dimension 5, in addition to the time variable.

In all examples considered in this section the operator $F(t,x,r,p,\gamma)$ does not depend on the $r-$variable. We shall then drop this variable from our notations, and we simply write the scheme as:
 \begin{equation} \label{schemenum}
 \begin{array}{l}
 v^h(T,.):=g
 ~~~\mbox{and}\\
 v^h(t_i,x):=\E[v^h(t_{i+1},\hat X_h^x)]
               + h F\left(t_i,x,\Dc_hv^h(t_i,x)\right)
 \end{array}
 \end{equation}
where
 \b*
\Dc_h\psi
 :=
 \left(\Dc_h^1\psi,\Dc_h^2\psi\right),
 \e*
and $\Dc_h^1$ and $\Dc_h^2$ are defined in Lemma \ref{lemintpart}. We recall from Remark \ref{remintpart} that:
 \begin{align} \label{schemes12}
 \Dc_{2h}^2 \varphi(t_i,x) 
 &\!\!=\!\!
 \E\!\left[\varphi(t_i\!+\!2h, \hat X^{t_i,x}_{2h})\!\!\left({\sigma^{\rm T}}\right)^{-1}\!\!
         \frac{(W_{t_i+h}\!-\!W_{t_i})(W_{t_i+h}\!-\!W_{t_i})^{\rm T}\!-\!h\mathbf{I}_d}{h^2}\;
       \!  \sigma^{-1}
   \right]
 \nonumber\\ 
 &\!\!=\!\!
 \E\!\left[  \Dc_h^1 \varphi(t_i\!+\!h,  \hat X^{t_i,x}_{h})\!\left({\sigma^{\rm T}}\right)^{-1}\!\Frac{W_{t_i+h}\!-\!W_{t_i}}{h}
   \right].
 \end{align}
The second representation is the one reported in \cite{cstv} where the present backward probabilistic scheme was first introduced. These two representations induce two different numerical schemes because once the expectation operator $\E$ is replaced by an approximation $\hat E^N$, equality does not hold anymore in the latter equation for finite $N$. In our numerical examples below, we provide results for both methods. The numerical schemes based on the first (resp. second) representation will be referred to as {\it scheme} 1 (resp. 2). An important outcome of our numerical experiments is that scheme 2 turns out to have a significantly better performance than scheme 1.

\begin{Remark}{\rm
The second scheme needs some final condition for $\Dc_h^1 \varphi(T,$ $ X^{T-h,x}_{h})$. Since $g$ is smooth in all our examples, we set this final condition to $\nabla g$. Since the second scheme turns out to have a better performnace, we may also use the final condition for $Z$ suggested by the first scheme.
}
\end{Remark}

We finally discuss the choice of the regression estimator in our implemented examples. Two methods have been used:
\begin{itemize}
\item The first method is the basis projection {\it a la} Longstaff and Schwartz \cite{longstaffschwartz}, as developed in \cite{glw}. We use regression functions with localized support : on each support the regression functions are chosen linear and the size of the support is adaptative according to the Monte Carlo distribution of the underlying process.
\item The second method is based on the Malliavin integration by parts formula as suggested in \cite{lionsr} and further developed in \cite{bet}. In particular, the optimal exponential localization function $\phi^k(y)=exp( -\eta^k y)$ in each direction $k$ is chosen as follows.
The optimal parameter $\eta_k$ is provided in \cite{bet} and should be chosen for each conditional expectation depending on $k$. Our numerical experiments however revealed that such optimal parameters do not provide sufficiently good performance, and more accurate results are obtained by choosing $\eta_k=5/\sqrt{\Delta t}$ for all values of $k$.
\end{itemize}

\subsection{Mean curvature flow problem}

The mean curvature flow equation describes the motion of a surface where each point moves along the inward normal direction with speed proportional to the mean curvature at that point. This geometric problem can be characterized as the zero-level set $S(t):=\{x\in\R^d: v(t,x)=0\}$ of a function $v(t,x)$ depending on time and space satisfying the geometric partial differential equation:
 \be\label{mcf}
 v_t - \Delta v + \frac{ Dv \cdot D^2 v Dv}{|Dv|^2} = 0
 &\mbox{and}&
 v(0,x) = g(x)
 \ee
and $g:\R^d\longrightarrow\R$ is a bounded Lipschitz-continuous function. We refer to \cite{sonertouzi} for more details on the mean curvature problem and the corresponding stochastic representation.

To model the motion of a sphere in $\R^d$ with radius $2R>0$, we take $g(x):=4R^2-|x|^2$ so that $g$ is positive inside the sphere and negative outside.
We first solve the sphere problem in dimension 3. In this case, it is well-known that the surface $S(t)$ is a sphere with a radius $R(t)=2\sqrt{ R^2 - t}$ for $t \in  (0,R^2)$.
Reversing time, we rewrite \reff{mcf} for $t\in(0,T)$ with $T = R^2$:
 \be\label{mcfrev}
 - v_t - \frac12\sigma^2 \Delta v + F(x,Dv,D^2v) = 0
 &\mbox{and}&
 v(T,x) = g(x),
 \ee

where
 \b*
 F(x,z,\gamma) &:=& \gamma\left(\frac12\sigma^2-1\right)+\frac{z\cdot\gamma z}{|z|^2}.
 \e*
We implement our Monte Carlo-finite differences scheme to provide an approximation $\hat v^h$ of the function $v$. As mentioned before, we implement four methods: Malliavin integration by parts-based or basis projection-based regression, and scheme 1 or 2 for the representation of the Hessian.

Given the approximation $\hat v^h$, we deduce an approximation of the surface $\hat S^h(t):=\{x\in\R^3: \hat v^h(t,x)=0)\}$ by using a dichotomic gradient descent method using the estimation of the gradient $\Dc^1 v$ estimated along the resolution. The dichotomy is stopped when the solution is localized within $0.01$ accuracy.

\begin{Remark}{\rm 
Of course the use of the gradient is not necessary in the present context where we know that $S(t)$ is a sphere at any time $t\in[0,T)$. The algorithm described above is designed to handle any type of geometry.
}
\end{Remark}

\begin{Remark}{\rm 
In our numerical experiments, the nonlinearity $F$ is truncated so that it is bounded by an arbitrary value taken equal to $200$.
}
\end{Remark}

Our numerical results show that Malliavin and basis projection methods give similar results. However, for a given number of sample paths, the basis projection method of \cite{glw} are slightly more accurate. Therefore, all results reported for this example correspond to the basis projection method.

Figure \ref{MCFsphere} provides results obtained with one million particles and $10 \times 10 \times 10$ mesh with a time step equal to $0.0125$. The diffusion coefficient $\sigma$ is taken to be either $1$ or $1.8$. We observe that results are better with $\sigma=1$.
We also observe that the error increases near time $0.25$ corresponding to an acceleration
of the dynamics of the phenomenon, and suggesting that a thinner time step should be used at the end of simulation.

\begin{figure}[h]
\centerline{
\includegraphics[width=10cm]{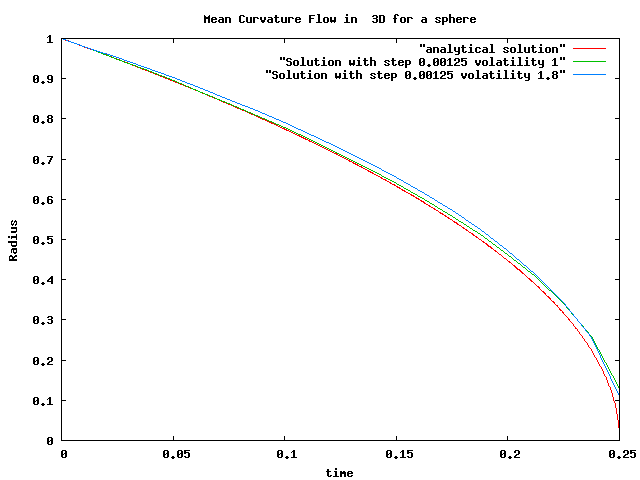}}
\caption{Solution of the mean curvature flow for the sphere problem}
\label{MCFsphere}
\end{figure}

Figure \ref{MCFspheremethod1} plots the difference between our calculation and the reference for scheme 1 and volatility 1 and 1.8 for varying time step.  The corresponding results with scheme 2 are reported in figure \ref{MCFspheremethod2}. We notice that some points at time $T=0.25$ are missing due to a non convergence of the gradient method for a diffusion $\sigma=1.8$. We observe that results for scheme 2 are slightly better than results for scheme 1. With $\sigma=1$, it takes 150 seconds on a Nehalem intel processor 2.9 GHz to obtain the result at time $t = 0.15$ with the regression method, while it takes  1500 seconds with the Malliavin method (notice that the dichotomy used with the gradient method is a very inefficient method).

\begin{figure}[h]
\centerline{
\includegraphics[width=6.5cm]{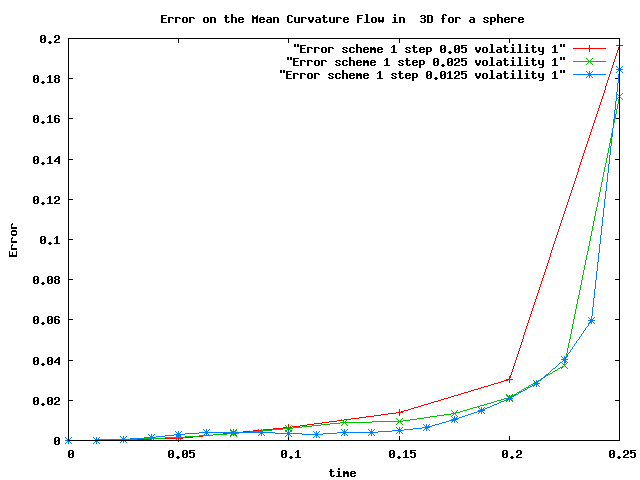}
\includegraphics[width=6.5cm]{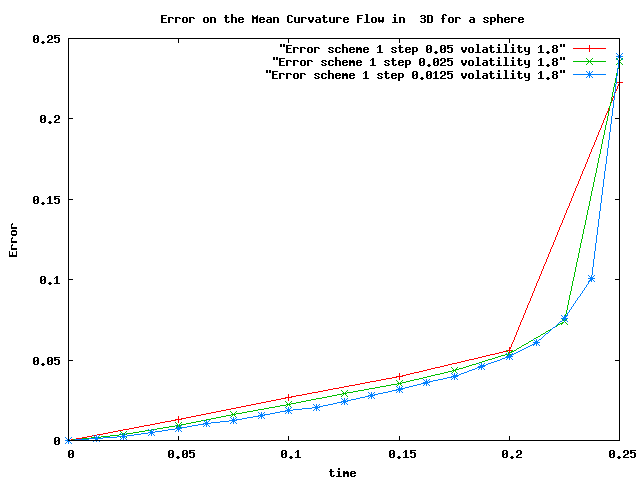}}
\caption{Mean curvature flow problem for different time step and diffusion: scheme 1}
\label{MCFspheremethod1}
\end{figure}

\begin{figure}[h]
\centerline{
\includegraphics[width=6.5cm]{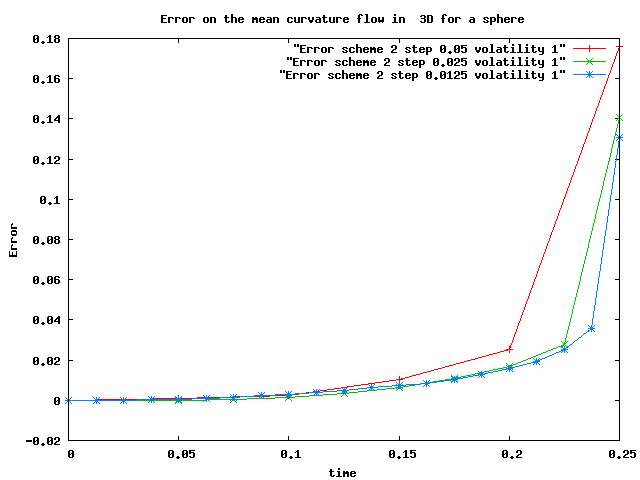}
\includegraphics[width=6.5cm]{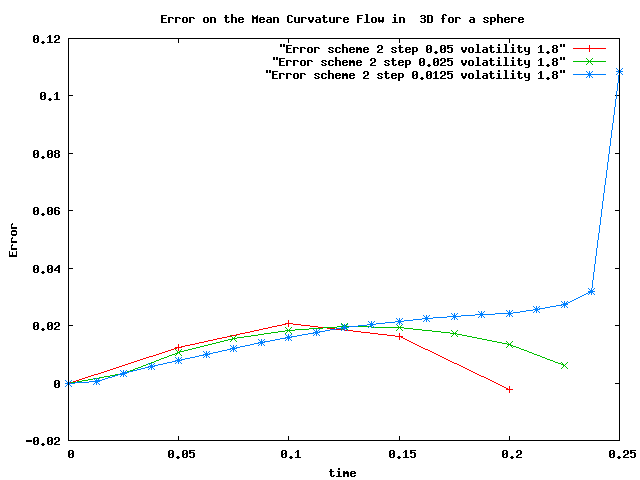}}
\caption{Mean curvature flow problem for different time step and diffusions: scheme 2}
\label{MCFspheremethod2}
\end{figure}

We finally report in Figure \ref{MCFtube} some numerical results for the mean curvature flow problem in dimension 2 with a more interesting geometry: the initial surface (i.e. the zero-level set for $v$) consists of two disks with unit radius, with centers positioned at -1.5 and 1.5 and connected by a stripe of unit width. We give the resulting deformation with scheme 2 for a diffusion $\sigma=1$, a time step $h=0.0125$, and one million particles. Once again, the Malliavin integration by parts based regression method and the basis projection method with
 $10 \times 10$ meshes produce similar results. We used $1024$ points to describe the surface.

One advantage of this method is the total parallelization that can be performed to solve 
the problem for different points on the surface : for the results given parallelization by Message Passing (MPI) was achieved.

\begin{figure}[h]
\centerline{
\includegraphics[width=10cm]{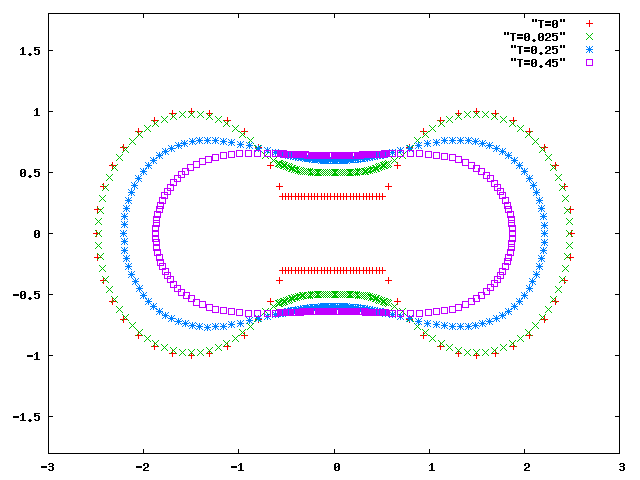}}
\caption{ Mean curvature flow problem in 2D}
\label{MCFtube}
\end{figure}

\subsection{Continuous-time portfolio optimization}

We next report an application to the continuous-time portfolio optimization problem in financial mathematics. Let $\{S_t,t\in[0,T]\}$ be an It\^o  process modeling the price evolution of $n$ financial securities. The investor chooses an adapted process $\{\theta_t,t\in[0,T]\}$ with values in $\R^n$, where $\theta_t^i$ is the amount invested in  the $i-$th security held at time $t$. In addition, the investor has access to a non-risky security (bank account) where the remaining part of his wealth is invested. The non-risky asset $S^0$ is defined by an adapted interest rates process $\{r_t,t\in[0,T]\}$, i.e. $dS^0_t=S^0_tr_tdt$, $t\in[0,1]$. Then, the dynamics of the wealth process is described by:
 $$
 dX^\theta_t
 =
 \theta_t\cdot \frac{dS_t}{S_t} +(X^\theta_t-\theta_t\cdot \1)\frac{d S^0_t}{S^0_t}
 = 
 \theta_t\cdot \frac{dS_t}{S_t}+(X^\theta_t-\theta_t\cdot \1)r_tdt,
 $$
where $\1=(1,\cdots,1)\in\R^d$.
Let $\Ac$ be the collection of all adapted processes $\theta$ with values in $\R^d$, which are integrable with respect to $S$ and such that the process $X^\theta$ is uniformly bounded from below. Given an absolute risk aversion coefficient $\eta>0$, the portfolio optimization problem is defined by:
 \be\label{prob-portefeuille}
 v_0
 &:=&
 \sup_{\theta\in\Ac} \E\left[-\exp\left(-\eta X^\theta_T\right)\right].
 \ee
Under fairly general conditions, this linear stochastic control problem can be characterized as the unique viscosity solution of the corresponding HJB equation. The main purpose of this subsection is to implement our Monte Carlo-finite differences scheme to derive an approximation of the solution of the fully nonlinear HJB equation in non-trivial situations where the state has a few dimensions. We shall first start by a two-dimensional example where an explicit solution of the problem is available. Then, we will present some results in a five dimensional situation.

\subsubsection{A two dimensional problem}

Let $d=1$, $r_t=0$ for all $t\in[0,1]$, and assume that the security price process is defined by the Heston model \cite{heston}:
 \b*
 dS_t &=& \mu S_tdt +  \sqrt{Y_t} S_t dW_t^{(1)}  
 \\
 dY_t &=&  k (m-Y_t)dt + c\sqrt{Y_t} \left(\rho dW_t^{(1)}+\sqrt{1-\rho^2}dW_t^{(2)}\right),
 \e*
where $W=(W^{(1)},W^{(2)})$ is a Brownian motion in $\R^2$.
In this context, it is easily seen that the portfolio optimization problem \reff{prob-portefeuille} does not depend on the state variable $s$. Given an initial state at the time origin $t$ given by $(X_t,Y_t)=(x,y)$, the value function $v(t,x,y)$ solves the HJB equation:
 \begin{equation}\label{hjb0}
 \begin{array}{rl}
 v(T,x,y) = - e^{-\eta x}
 ~\mbox{and}~
 0 
 =&\!\!\!\!
 - v_t -  k (m-y) v_y - \frac{1}{2}c^2 y v_{yy}\\
&\hspace*{1.7cm}
 -\Sup_{\theta\in\R} \Bigl(\frac12\theta^2yv_{xx}+\theta(\mu v_x+\rho cy v_{xy})\Bigr) 
 \\
 =&\!\!\!\!
 - v_t -  k (m-y) v_y - \frac{1}{2}c^2 y v_{yy}
 + \Frac{(\mu v_x + \rho cy v_{xy})^2} 
        {2 y v_{xx}}.
 \end{array}
 \end{equation}
A quasi explicit solution of this problem was provided by Zariphopoulou \cite{zari}:
 \be\label{zariphopoulou}
 v(t,x,y)=-e^{-\eta x} \left\| \exp\left(-\frac{1}{2}\int_t^T \frac{\mu^2}{\tilde Y_s}ds\right)
                       \right\|_{\L^{1-\rho^2}}
 \ee
where the process $\tilde Y$ is defined by
 \b*
 \tilde Y_t=y
 &\mbox{and}&
 d\tilde Y_t = ( k (m-\tilde Y_t)- \mu c\rho)dt +  c\sqrt{\tilde Y_t} dW_t.
 \e*
In order to implement our Monte Carlo-finite differences scheme, we re-write \reff{hjb0} as:
 \begin{equation}\label{hjb1}
 - v_t -  k (m-y) v_y - \frac{1}{2}c^2 yv_{yy} - \frac{1}{2}\sigma^2 v_{xx} 
 + F\left(y,Dv,D^2v\right) = 0, 
 ~
 v(T,x,y) = - e^{-\eta x},
 \end{equation}
where $\sigma>0$ and the nonlinearity $F:\R\times\R^2\times\S_2$ is given by:
 \b*
 F(y,z,\gamma) 
 &=& 
 \frac{1}{2} \sigma^2 \gamma_{11} +\frac{(\mu z_1 + \rho c y \gamma_{12})^2}
                                        {2 y \gamma_{11}}.
 \e*
Notice that the nonlinearity $F$ does not to satisfy Assumption {\bf F}, we consider the truncated nonlinearity:
 \b*
 F_{\eps,M}(y,z,\gamma) 
\!\! &\!\!:=\!\!&\!\!
 \frac{1}{2} \sigma^2 \gamma_{11} 
 -
 \sup_{\eps\le\theta\le M} \left(\frac12\theta^2(y\vee\eps)\gamma_{11}+\theta(\mu z_1+\rho c(y\vee\eps) \gamma_{12}\right),
 \e*
for some $\eps, n>0$ jointly chosen with $\sigma$ so that Assumption {\bf F} holds true.
Under this form, the forward two-dimensional diffusion is defined by:
 \be\label{sde}
 dX^{(1)}_t =  \sigma dW^{(1)}_t, 
 &\mbox{and}&
 dX^{(2)}_t =   k (m-X^{(2)}_t) dt  +  c\sqrt{X^{(2)}_t} dW^{(2)}_t.
 \ee
In order to guarantee the non-negativity of the discrete-time approximation of the process $X^{(2)}$, we use the implicit Milstein scheme \cite{kahl}:
 \be\label{milstein}
 X_{n}^{(2)} 
 &=& 
 \frac{X_{n-1}^{(2)} +  k   m \Delta t + c  \sqrt{X_{n-1}^{(2)}} \xi_{n} \sqrt{ \Delta t} 
                + \frac{1}{4} c ^2  \Delta (\xi_n^2-1)}{1 +  k  \Delta t}
 \ee
where $(\xi_n)_{n\ge 1}$ is a sequence of independent random variable with distribution $\mathbf{N}(0,1)$.

Our numerical results correspond to the following values of the parameter: $\mu =0.15$, $c  = 0.2$, $ k  = 0.1$, $m  = 0.3$, $Y_0 =  m$, $\rho =0$. The initial value of the portfolio is $x_0=1$, the maturity $T$ is taken equal to one year. With this parameters, the value function is computed from the quasi-explicit formula \reff{zariphopoulou} to be $v_0=-0.3534$.

We also choose $M=40$ for the truncation of the nonlinearity. This choice turned out to be critical as an initial choice of $M=10$ produced an important bias in the results.

The two schemes have been tested with the Malliavin and basis projection methods. The latter was applied with $40 \times 10$ basis functions. We provide numerical results corresponding to 2 millions particles. 
Our numerical results show that the Malliavin and the basis projection methods produce very similar results, and achieve a good accuracy: with 2  millions particles, we calculate the variance of our estimates by performing 100 independent calculations:
\begin{itemize}
\item the results of the Malliavin method exhibit a standard deviation smaller than $0.005$ for scheme one (except for a step equal to  $0.025$ and  a volatility equal to $1.2$ where standard deviation jumped to $0.038$), $0.002$ for scheme two  with a computing time of 378 seconds for 40 time steps, 
\item the results of the basis projection method exhibit a standard deviation smaller than $0.002$ for scheme 1 and $0.0009$ for scheme two with a computing time of 114 seconds for 40 time steps. 
 \end{itemize}

Figure \ref{figfin1} provides the plots of the errors obtained by the integration by parts-based regression with Schemes one and two. 
All solutions have been calculated as the average of 100 calculations. We first observe that for a small diffusion coefficient $\sigma=0.2$, the numerical performance of the algorithm is very poor: surprisingly, the error increases as the time step shrinks to zero and the method seems to be biased. This numerical result hints that the requirement that the diffusion should dominate the nonlinearity in Theorem \ref{thmconv}, might be a sharp condition.
\begin{figure}[h]
\centerline{
\includegraphics[width=6.5cm]{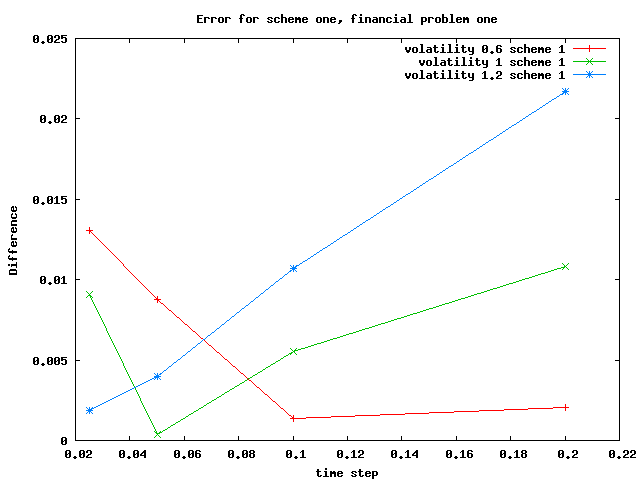}
\includegraphics[width=6.5cm]{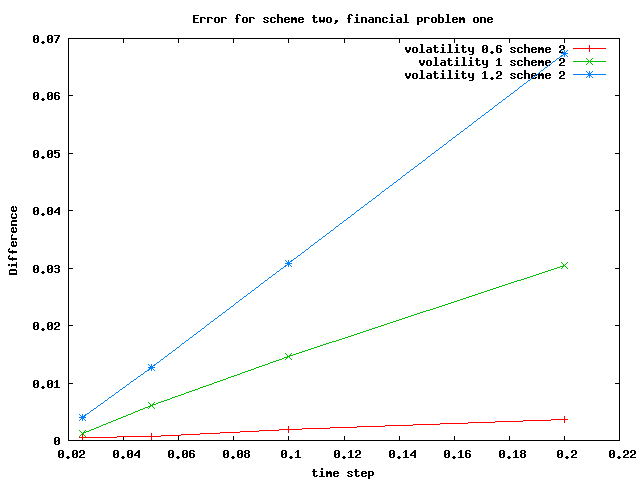}}
\caption{ Difference between calculation and reference for scheme one and two }
\label{figfin1}
\end{figure}
We also observe that scheme one has a persistent bias even for a very small time step, while scheme two exhibits a better convergence towards the solution.

\subsubsection{A five dimensional example}

We now let $n=2$, and we assume that the interest rate process is defined by the Ornstein-Uhlenbeck process:
 \b*
 dr_t &=& \kappa(b-r_t)dt+\zeta dW^{(0)}_t.
 \e*
While the price process of the second security is defined by a Heston model, the first security's price process is defined by a CEV-SV models, see e.g. \cite{lord} for a presentation of these models and their simulation:
 \b*
 dS^{(i)}_t
 &=& 
 \mu_i S^{(i)}_tdt + \sigma_{i} \sqrt{Y^{(i)}_t} {S^{(i)}_t}^{\beta_i} dW^{(i,1)}_t,~~\beta_2=1,
 \\
 dY^{(i)}_t 
 &=&  k_i\left(m_i -Y^{(i)}_t\right)dt +  c_i  \sqrt{Y^{(i)}_t} dW^{(i,2)}_t 
 \e*
where $\left(W^{(0)},W^{(1,1)},W^{(1,2)},W^{(2,1)},W^{(2,2)}\right)$ is a Brownian motion in $\R^5$, and for simplicity we considered a zero-correlation between the security price process and its volatility process. 

Since $\beta_2=1$, the value function of the portfolio optimization problem \reff{prob-portefeuille} does not depend on the $s^{(2)}-$variable. 
Given an initial state $(X_t,r_t,S^{(1)}_t$ $,Y^{(1)}_t,Y^{(2)}_t)=(x,r,s_1,y_1,y_2)$ at the time origin $t$, the value function $v(t,x,r,$ $s_1,y_1,y_2)$ satisfies the HJB equation:
 \begin{align}
 0&=-v_t - (\mathbf{L}^r +\mathbf{L}^Y +\mathbf{L}^{S^1})v - r x v_x 
 \nonumber\\
  & -\sup_{\theta_1,\theta_2} 
      \left\{\theta_1\cdot(\mu\!-\!r\1)v_x
             \!+\!\theta_1 \sigma_1^2 y_1 s_1^{2 \beta_1-1} v_{xs_1}
             \!+\!\frac12(\theta_1^2\sigma_1^2 y_1 s_1^{2\beta_1-2}
             \!+\!\theta_2^2\sigma_2^2y_2)v_{xx}
      \right\}
 \nonumber\\
  &=
  -v_t - (\mathbf{L}^r +\mathbf{L}^Y+\mathbf{L}^{S^1})v 
  - rx v_x 
  \nonumber\\
  &  +\Frac{((\mu_1- r)v_x+\sigma_1^2 y_1 s_1^{2 \beta_1-1} v_{xs_1})^2}
                 {2 \sigma_1^2 y_1 s_1^{2\beta_1-2}v_{xx}}
  +\Frac{((\mu_2-r)v_x)^2}{2 \sigma_2^2 y_2 v_{xx}}
 \label{equation5D}
 \end{align}
where
 \b*
 \mathbf{L}^r v=\kappa(b-r)v_r+\frac12\zeta^2 v_{rr},
 &&
 \mathbf{L}^Y v=\sum_{i=1}^2k_i\left(m_i -y_i\right)v_{y_i} +  \frac12 c_i^2  y_iv_{y_i y_i},
 \e*
 \b*
 \mbox{and}~~\mathbf{L}^{S^1}v
 &=&
 \mu_1 s_1v_{s_1}-\frac12\sigma_1^2 s_1 y_1 v_{s_1s_1}.
 \e*
In order to implement our Monte Carlo-finite differences scheme, we re-write \reff{equation5D} as:
 \begin{equation}\label{hjb2}
\begin{array}{l}
 - v_t - (\mathbf{L}^r +\mathbf{L}^Y +\mathbf{L}^{S^1})v - \frac{1}{2}\sigma^2 v_{xx} 
 + F\left((x,r,s_1,y_1,y_2),Dv,D^2v\right) = 0, \\
 ~~
 v(T,x,r,s_1,y_1,y_2) = - e^{-\eta x},
\end{array}
 \end{equation}
where $\sigma>0$, and the nonlinearity $F:\R^5\times\R^5\times\S_2$ is given by:
 \b*
 F(u,z,\gamma) 
 \!\!&\!\!\!=\!\!\!& \!\!
 \frac{1}{2} \sigma^2 \gamma_{11} \! -\! x_1 x_2 z_1\! +\!
 \Frac{((\mu_1\!-\! x_2) z_1\!+\!\sigma_1^2 x_4 x_3^{2 \beta_1-1} \gamma_{1,3})^2}{2 \sigma_1^2 x_4 x_3^{2\beta_1-2} \gamma_{11}}
 \!+\!\Frac{((\mu_2\!-\!x_2) z_1)^2}{2 \sigma_2^2 x_5 \gamma_{11}},
 \e*
where $u=(x_1,\cdots,x_5)$.
We next consider the truncated nonlinearity:
 \b*
 F_{\eps,M}(u,z,\gamma)
 \!\!&\!\!\!:=\!\!\!&\!\!
 \frac{1}{2} \sigma^2 \gamma_{11}  \!- \!x_1 x_2 z_1 \!+\!
 \!\!\!\!\!\sup_{\eps\le|\theta|\le M} \!
      \Big\{(\theta\cdot(\mu\!-\!r\1)z_1
   \!+\! \theta_2^2\sigma_2^2(x_5\vee\eps))\gamma_{11}\\
\!\!&\!\!\!\!&             \!+\!\theta_1 \sigma_1^2 (x_4\vee\eps) (x_3\vee\eps)^{2 \beta_1-1} \gamma_{13}
   \!+\!\frac12(\theta_1^2\sigma_1^2 (x_3\vee\eps) (x_4\vee\eps)^{2\beta_1-2}
                \Big\},          
 \e*
where $\eps,M>0$ are jointly chosen with $\sigma$ so that Assumption {\bf F} holds true. Under this form, the forward two-dimensional diffusion is defined by:
 \begin{equation}\label{sde1}
\begin{array}{lll}
 dX^{(1)}_t =  \sigma dW^{(0)}_t,
 \\
 dX^{(2)}_t =   \kappa(b-X^{(2)}_t)dt+\zeta dW^{(1)}_t, \\
  dX^{(3)}_t =    \mu_1 X^{(3)}_t dt + \sigma_{1} \sqrt{X^{(4)}_t} {X^{(3)}_t}^{\beta_1} dW^{(1,1)}_t,
 \\
  dX^{(4)}_t =   k_1 (m_1-X^{(4)}_t) dt  +  c_1\sqrt{X^{(4)}_t} dW^{(1,2)}_t, \\
  dX^{(5)}_t =   k_2 (m_2-X^{(5)}_t) dt  +  c_2\sqrt{X^{(5)}_t} dW^{(2,2)}_t. 
 & &
\end{array}
\end{equation}
The component $X^{(2)}_t$ is simulated according to the exact discretization:
\b*
X^{(2)}_{t_n}
&=&
b+e^{-k\Delta t}\left(X^{(2)}_{t_{n-1}}-b\right)
+\zeta\sqrt{\frac{1- \exp(- 2 \kappa \Delta t)}{2 \kappa}}\xi_n,
\e*
where $(\xi_n)_{n\ge 1}$ is a sequence of independent random variable with distribution $\mathbf{N}(0,1)$.
The following scheme for the price of the asset guarantees non-negativity (see \cite{andersen}) :

\b*
\ln X_n^{(3)}\!\! =\!\!  \ln X_{n-1}^{(3)} \!\!+\!\!\left (\!\!\mu_1 \!-\! \frac{1}{2} \sigma^2_{1} \!\left(\!X_{n-1}^{(3)}\!\right)^{\!\!2(\beta_1\!-\!1)} \!\!\!\! X_{n-1}^{(4)} \!\!\right)\!\!\Delta t \!+\! 
  \sigma_{1} \!\! \left(X_{n-1}^{(3)}\right)^{\!\!\beta_i-1} \!\!\!\!\!\sqrt{\!\!X_{n-1}^{(4)}} \Delta W_n^{(1,2)}
\e*
where $\Delta W_n^{(1,2)}:=W_n^{(1,2)} - W_{n-1}^{(1,2)}$. We take the following parameters  $\mu_1 = 0.10$, $\sigma_1 =0.3$, $\beta_1 =0.5$ for the first asset,
$ k _1 =0.1$, $m _1 = 1.$, $c_1 = 0.1$ for the diffusion process of the first asset. The second asset is defined by the same parameters as in the two dimensional example: $\mu_2=0.15$, $c_2=0.2$, $m=0.3$ and $Y^{(2)}_0=m$. As for the interest rate model we take $b=0.07$, $X_0^{(2)} = b$, $\zeta= 0.3$.

The initial values of the portfolio the assets prices are all set to 1. For this test case we first use the basis projection regression method with $4 \times 4 \times 4 \times 4 \times 10$ meshes and three millions particles which, for example, takes 520 seconds for 20 time steps.
Figure \ref{figfin2} contains the plot of the solution obtained by scheme 2, with different time steps. We only provide results for the implementation of scheme 1 with a coarse time step, because the method was diverging with a thinner time step. We observe that there is still a difference for very thin time step with the three considered values of the diffusion. This seems to indicate that more particles and  more meshes are needed. While doing many calculation we observed that for the thinner time step mesh, the solution sometimes diverges. 
We therefore report the results corresponding to thirty millions particles with $4 \times 4 \times 4 \times 4 \times 40$ meshes. 
First we notice that with this discretization all results are converging as time step goes to zero: the exact solution seems to be very closed to $- 0.258$. During our experiments with thirty millions particles, the scheme was always converging with a very low variance on the results. A single calculation takes 5100 seconds with 20 time steps.

\begin{Remark}{\rm
With thirty millions particles, the memory needed forced us to use 64-bit processors with more than four gigabytes of memory.
}
\end{Remark}

\begin{figure}[h]
\centerline{
\includegraphics[width=6.5cm]{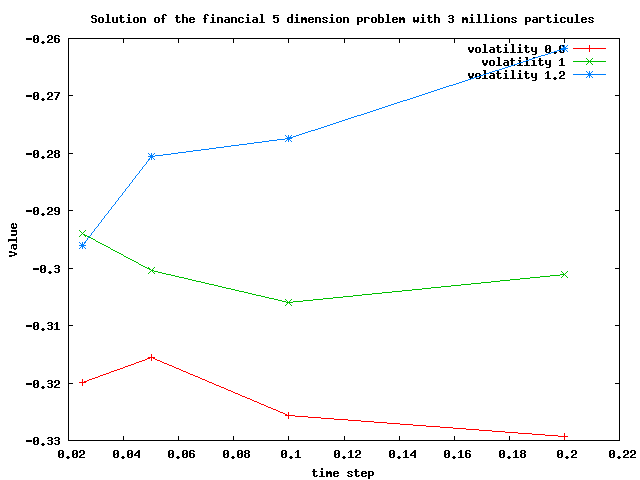}
\includegraphics[width=6.5cm]{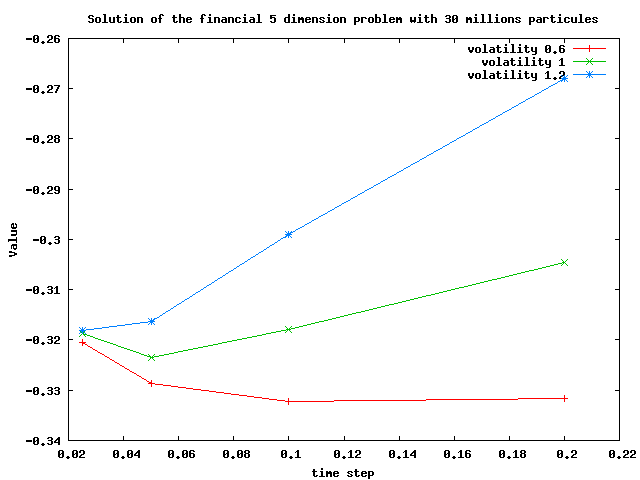}
}
\caption{Five dimensional financial problem and its results for different volatilities with 3 millions and 30 millions particles }
\label{figfin2}
\end{figure}

\subsubsection{Conclusion on numerical results}

The Monte Carlo-Finite Differences algorithm has been implemented with both schemes suggested by \reff{schemes12}, using the basis projection and Malliavin regression methods. Our numerical experiments reveal that the second scheme performs better both in term of results and time of calculation for a given number of particles, independently of the regeression method. 

We also provided numerical results for different choices of the diffusion parameter in the Monte Carlo step. We observed that small diffusion coefficients lead to poor results, which hints that the condition that the diffusion must dominate the nonlinearity in Assumption {\bf F} (iii) may be sharp. On the other hand, we also observed that large diffusions require a high refinement of the meshes meshes, and large number of particles, leading to a high computational time.

Finally, let us notice that a reasonable choice of the diffusion could be time and state dependent, as in the classical importance sampling method. We have not tried any experiment in this direction, and we hope to have some theoretical results on how to choose optimally the drift and the diffusion coefficient of the Monte Carlo step.

\end{document}